\documentstyle[12pt]{amsart}
\begin{document}

\def\CC{\bf C}
\def\RR{\bf R}
\def\C{\hbox{C}}
\def\g{\gamma}
\def\rh{\rho}
\def\DD{\D_{\a^*}}
\def\ale{\mathrel{\mathop{<}\limits_{\sim}}}
\def\ld{\lambda}
\def\sm {\setminus}
\def\O{\Omega}
\def\b{\beta}
\def\o{\omega}
\def\D{\Delta}
\def\ss{\subset}
\def\1{{(0,\gamma)\OO \|\gamma\|}}
\def\L{\Lambda}
\def\ic{e^{i\theta}}
\def\pd{\phi(\Delta)}
\def\h{\hbox}
\def\df{{d\Phi^*\OO d\xi}}
\def\q{Q(X,\xi)}     
\def\wt{\widetilde}
\def\ra{\rightarrow}
\def\ih{I_{2n+1}+i{\partial H(X)\over\partial X }}
\def\PP{\partial}
\def\OO{\over}
\def\d{\delta}
\def\l{\lambda}
\def\Sqt{\sqrt}
\def\endpf${$\qed}
\def\H{\hat}
\def\a{\alpha}
\def\d{\delta}
\def\V{\h{Val}}
\def\Ol{\overline}
\def\OO{\over}
\def\R{\h{Re}}
\def\I{\h{Im}}
\def\C{\hbox{C}}
\hsize=6.0truein
\vsize=8.0truein
\def\D{\Delta}
\def\d{\delta}
\def\a{\alpha}
\def\zn{z_{n+1}}
\def\Sqt{\sqrt}

\title[Regularity of Holomorphic Correspondences] 
{Regularity of Holomorphic Correspondences and
Applications to the Mapping Problems}

\author{Xiaojun Huang}
\address{
Department of Mathematics, The University of Chicago, Chicago, IL 60637, USA}
\curraddr{
Mathematical Sciences Research Institute, 1000 Centennial Drive, Berkeley, CA 94720, USA} 
\email{xhuang@@math.uchicago.edu; xhuang@@msri.org}

\thanks{Supported in part by NSF DMS-9500881.
Research at MSRI is supported in part by NSF DMS-9022140.}

\begin{abstract}
We study the regularity results of holomorphic
 correspondences. As an application, we combine it with certain
  recently developed  methods to obtain
the extension theorem for  proper holomorphic mappings between domains
with real analytic boundaries in the complex 2-space.
\end{abstract}

\maketitle

\section*{Introduction} 

This paper is
concerned with the boundary regularity problem for holomorphic
 mappings 
between  domains in complex euclidean spaces. Our main
purpose is to present certain results  related to the following
mapping conjecture in several complex variables:
\bigbreak
{\bf Conjecture 0.1}: Let $M_1$ and $M_2$ 
be two real analytic hypersurfaces of finite D-type in $\bf C^n$ with $n>1$.
Suppose that $D\subset {\bf C^n}$ is a domain with $M_1$ as part
of its real analytic boundary.
 Also, suppose that $f$ is a holomorphic mapping
from $D$ into $\bf C^n$, that is continuous up to $D\cup M_1$ and
maps $M_1$ into $M_2$.
Then
 $f$ admits a holomorphic extension across $M_1$.
\bigbreak
 Here, we recall that a real
analytic hypersurface is  of finite D-type [Da] 
if and only if it does not contain
any non-trivial complex analytic variety.

The extensive study of such  kinds of problems
 was initiated by the foundational work
of Fefferman in the case of biholomorphic mappings
between strongly pseudoconvex domains.
Recent related papers
 include those  by
Lewy [Le], Pinchuk [Pi1], Webster [We], Diederich- Webster [DW], Bell [Be1] 
[Be3],
Baouendi-Jacobowitz-Treves [BJT], Bedford-Bell [BB], 
Baouendi - Bell- Rothschild [BBR], Baouendi- Rothschild [BR1] [BR2] [BR3] [BR4], 
Bell-Catlin [BC1] [BC2], Diederich-Fornaess [DF1] [DF2] [DF3] [DF4],
 Pinchuk-Tsyganov
[PT], 
Diederich- Fornaess- Ye [DFY],
 Baouendi- Huang- Rothschild [BHR1] [BHR2], Huang -  Pan [HP], Pan [Pan],
 and the references therein.
For  detailed surveys of this investigation, 
we refer the reader to the articles by Forstneric [For], 
Bedford [Bed], and Bell-Narashimhan [BN].

When $M_1$ and $M_2$ are pseudoconvex 
and when $f$ is proper, a  powerful method 
is  to use   Condition R
introduced by Bell-Ligocka (see [Be2], for example) and the solutions 
to the $\overline {\PP}$-Neumann problems 
to achieve  the smooth extension. Then one  applies 
various  versions of the reflection principle 
to obtain the holomorphic extension ([BJT] [BBR] [BR1] [DF1]). However,
because of an example by  Barrett [Ba1], it is known 
that 
a  different approach needs to be employed 
when the assumption of pseudoconvexity  is dropped. 

In [We], Webster used  a class of invariant varieties,   
introduced by  Segre  [Se], to analyze biholomorphisms
between non-degenerate algebraic hypersurfaces.
Later, it became clear  that these invariant objects called Segre surfaces
 may also
be useful to  study 
the  above problem even in the non-pseudoconvex
case (for example, see  [DW], [DF2], [DFY]).

The plausible method based of the use of Segre surfaces
might  be carried out in two steps.
First one proves the holomorphic extension under the assumption that
the map extends as a
holomorphic correspondence (see the following section for
a  precise definition of this terminology). Then, one
reduces the general situation to the above case.
A   result along these lines  was obtained
in [DF2], where it was shown that a biholomorphism between two
real analytic domains  
in $\bf C^2$, that extends as a correspondence, 
admits a holomorphic extension up to the boundary. 
(See also related work in [BB], [Art], [BBR], [DF1] when the map is smooth). 
This, in particular, implies that
any biholomorphic map between two  
bounded algebraic domains in $\bf C^2$
extends holomorphically across the boundary [DF2]. 
The idea used in [DF2]  
depends strongly on a theorem of Barrett [Ba2], and does not seem 
 adaptable to the 
proper mapping case.

In this paper, we will establish the smoothness results  for
holomorphic correspondences coming from general holomorphic mappings
(Theorem 1.1$'$).
 As an immediate application,
we obtain a solution of Conjecture 0.1 in case $n=2$, $M_1$ and $M_2$ are
alegbraic, and $f$ is proper.
One of the other main results
indicates 
that any continuous CR mapping between two variety-free
real algebraic hypersurfaces in $\bf C^2$ is real analytic.
  This, in particular, gives a complete
solution to Conjecture 0.1 in the 2-dimensional   algebraic case.
Previously, Conjecture 0.1 was only known to hold in case hypersurfaces are
strongly pseudoconvex.  Our result seems
to give
 the first Schwarz reflection principle for {\it continuous CR mappings}
between a large class of hypersurfaces which may not be pseudoconvex,
though much more has been well understood when the map is smooth and proper
(see the papers mentioned above). 
Indeed,  it is a wide open question to answer
 whether Conjecture 1.1 holds when
the hypersurfaces are pseudoconvex  of finite type ([BC]), or
  when $f$ is smooth up to $D\cup M_1$ but
the hypersurfaces may not be pseudoconvex ([BR2]).

Recently,   based on the aforementioned result in [DF2] 
for biholomorphisms and the machinery developed in [We],
Diederich -Pinchuk ( [DP]) announced a proof 
 that any biholomorphic map between
bounded real analytic  domains in $\bf C^2$ extends  
as a holomorphic correspondence across the boundary (see also closely
related work in
[DFY]).
This together with the result in [DF2] implies the holomorphic
extension for biholomorphisms between real analytic domains in $\bf C^2$.
With the main  
result in the present paper  as our major tool (Theorem 1.1 or
Theorem 1.1$'$), and  with the help
 of 
ideas the author learned from [DFY] and [DP], 
we will  obtain, in the the last section of 
this paper, a solution of Conjecture 0.1 in case $n=2$ and
$f$ is proper. Indeed, once Theorem 1.1 is known, the rest of the argument
of the proof of Theorem 1.3 is a  modification of the known methods and ideas
for biholomorphic maps as appeared in  [DFY], and, in particular, [DP].
However, since [DP] may not be  presently available,
 we give a certain detailed discussion
on this matter ($\S 7$)  for completeness of the proof of Theorem 1.3.

\bigbreak
Our method of the proof of Theorem 1.1,
 which is different from that in [DF2], can be described as 
follows: We first show that an attached reflection function
(see, in particular, [BJT], [BBR])
extends holomorphically across the hypersurface. 
 Here the main
idea is to use the weak version of the edge of the wedge theorem
(see Proposition 7.1 of [BHR1] or [Hu]).
Using this reflection function,
 we will  connect the branching locus of the map with the branching locus
of the Segre varieties of the target hypersurface. 
 This then gives  us  many restrictions
 on the branching behaviors of the map. Finally the above mentioned
results and several  technical
lemmas  (in particular, a modified version of the
Baouendi-Rothschild Hopf lemma and a preservation principle) give the proof of Theorem 1.1.
The proof of Theorem 1.2 is  an easy corollary of Theorem 1.1 with
the known results. Once one has Theorem 1.1 as the major tool, the proof of 
Theorem 1.3 can be achieved with the help of 
 the existing ideas (see [DF1] [DFY]
and, in particular, [DP]). 

We would like to mention that
certain notations in this paper (especially,  those in $\S 7$) are adapted from
[DFY] and, in particular, [DP].

\bigbreak

{\bf Acknowledgment}: The author is grateful to 
Professors S. Webster and R. Nara\-shim\-han  for their many stimulating
conversations
(in particular, certain discussions related to Lemma 5.3), and
constant encouragement. He  would like
to thank Professors S. Baouendi,  S. Krantz, 
and L. Rothschild for
their constant help and interest in this work. The author also thanks Professor
S. Pinchuk for his   explanation of his work with K. Diederich [DP]
and for showing his lecture notes [DP] (on the construction of
holomorphic correspondences for biholomorphic mappings). 
\bigbreak
{\it Added in revision}:
This is the revised version of the author's previous paper with the same title,
in which Theorem 1.1 was proved only in case
 $n=1$ or $M_1$ and $M_2$ are decoupled.
 Results related to Theorem 1.1$'$ were also announced by Diederich-Pinchuk
in  the workshop on Analytic and PDE methods in SCVs at MSRI (November, 1995).

\section{Statement of main results and related observations} 

To state our main results, we need the following preliminary notation
(see  [BB], [DF2], and, in particular, [DFY], [DP]):

Let $M_1$ and $M_2$ be two real analytic hypersurfaces of finite D-type
in $\bf C^{n+1}$ ($n\ge 1$), and let $p$ be a point in $M_1$.
A  CR mapping $f$ from $M_1$ into $M_2$ is said to extend as 
 a {\it holomorphic correspondence}
to a neighborhood of $p$ 
if (a) $f$ is continuous near $p$; (b)
for each component $f_j$ of $f$, there is a 
polynomial ${\cal P}_j(z;X)$
in $X$ with leading coefficient 1 and other coefficients 
holomorphic in  a small neighborhood $U_1$ of $p$,
 such that ${\cal P}_j(z;f_j(z))
= 0$ for $z\in U_1\cap M_1$. 
It can be shown that $f$ extends as a holomorphic correspondence at $p$
if and only if there are a small neighborhood $U_1$ of $p$, a neighborhood
$U_2$ of $M_2$, and an analytic sub-variety $V$ of $U_1\times U_2$
such that 
(a) $ V\supset \Gamma_f=\{(z,w): w=f(z),\ z\in M_1\cap U_1\}$ 
(b) the natural projection of $V$ to  $U_1$ 
is locally proper.

We also recall that
a real analytic hypersurface $M\subset {\bf C^{n+1}}$ is called {\it rigid} if,
after a holomorphic change of coordinates, it can
be defined by a real analytic function of 
the form $\rho(z,\overline {z})=\zn+\overline {\zn}+\chi(z',\overline {z'})$, where
$z'=(z_1,\cdots,z_n)$.

\bigbreak
{\bf Theorem 1.1}:
Let $M_1$ and $M_2$ be two real analytic hypersurfaces 
of finite D-type  in $\bf C^{n+1}$ ($n\ge 1$). Let $p\in M_1$.
Suppose that $f$ is a  CR mapping from $M_1$ into $M_2$ that extends as 
a holomorphic
correspondence to a neighborhood of $p$. Then $f$ admits a holomorphic
extension near $p$, if either $n=1$ or $M_1$ and $M_2$ are rigid.
\bigbreak
 See also Theorem 1.1$'$ of Section 5 for
a more general version of this result and notice that
 in whole paper, only in $\S 5$, 
what we call Condition S is used.

Recall that a smooth hypersurface is called algebraic, if it can be
 defined by a real polynomial.
Using certain known results in [We] and [DF2],  one can apply Theorem 1.1 to
obtain the following:
\bigbreak
{\bf Theorem 1.2}: Let $M_1$ and $M_2$ be two algebraic hypersurfaces of finite
type in $\bf C^2$. Let $f$ be a continuous 
CR mapping from $M_1$ to $M_2$. Then
$f$ is real analytic on $M_1$.

\bigbreak
{\bf Corollary 1.2$'$}: Let $D$ and $D'$ be two smoothly bounded
 domains with real algebraic boundaries
in $\bf C^2$. Let $f$ be a proper holomorphic
mapping from $D$ to $D'$. Then $f$ admits a holomorphic extension
to a neighborhood of $\overline {D}$.

\bigbreak
We remark that Theorem 1.2 is sharp in a certain sense. 
In fact,   using the inner functions
defined on the ball in $\bf C^2$, one can construct
a bounded CR mapping between  spheres in $\bf C^2$ 
which is even not continuous
at any point. 
 It is worthwhile to mention that
the map in Theorem 1.2
 is not assumed to be the boundary value of some globally defined 
proper holomorphic
map. As is known, it is always a difficult problem to understand 
when a CR mapping
comes from a globally defined proper holomorphic map. 
Also, we notice that 
 when the map in Theorem 1.2
 is assumed to be in the smoothness class $C^k$ with
 $k>>1$,  then much less restriction to the hypersurfaces is required
to obtain the real analyticity ( [BHR1]).

Using the two dimensional homeomorphic version of 
Theorem 1.1 (which has been established in [DF2]),
Diederich-Pinchuk recently announced that any biholomorphic map between real
analytic bounded domains in $\bf C^2$ admits a correspondence extension
across the boundary [DP] (see also
the related  work in [DFY]).  
 With Theorem 1.1 at our disposal and with the help of an idea 
of Diederich-Pinchuk [DP] on the application of the Bishop extension lemma,
 we will prove in $\S 7$
 the following result  (see [BBR] for the
 pseudoconvex case, and Theorem 1.3$''$ in $\S 7$ for a local version of this
result).

\bigbreak
{\bf Theorem 1.3}:
Let $D$ and $D'$ be two bounded  domains 
in $\bf C^2$ with real analytic boundaries.
 Suppose that $f$ is a proper holomorphic map from $D$ to
$D'$, which extends continuously across $\PP D$. Then $f$ extends 
holomorphically
across $\PP D$.
\bigbreak
 In higher dimensions, we will be content in this paper with the
following application of Theorem 1.1 (see Theorem 1.1 $'$ for a more general
result on this):

 \bigbreak
 {\bf Theorem 1.4}: Let $D$ and $D'$ be two smoothly bounded
 algebraic domains in $\bf C^n$ and let $f$ be a proper holomorphic map from
 $D$ to $D'$. If $D$ and $D'$ are rigid 
then $f$ admits a holomorphic extension to a neighborhood
of $\overline {D}$.

\bigbreak
The organization of the remainder of this paper
  is as follows:  $\S 2$ through $\S 5$
is devoted to the proof of Theorem 1.1. 
In $\S 2$, we will study
the holomorphic extendibility of a $\l$-function associated 
with the map and the target
surface. This function  was  used in  various different forms in 
[Le], [Pi], [BJT], [BBR], and [BR1]
to study  smooth CR-mappings (see, in particular, the work in [BBR] and [BR1]).
In $\S 3$, we will present the connections
between the branching properties of the map and the invariant varieties 
(Segre-surfaces). In $\S 4$, we will prove a localized
Hopf lemma by 
modifying an argument that 
appeared in [BR2], which will then be  used to prove
a    preservation principle for our (multiple-valued) map. 
Section 5 is devoted to the
proof of Theorem 1.1,  using results  established in the previous
sections. In $\S 6$, we will complete the proofs of 
Theorem 1.2 and Theorem 1.4. With 
Theorem 1.1 at our disposal and using some known
 ideas (in particular, those appeared in [DP]),
 in the last
section ($\S 7$), we will present a proof of Theorem 1.3.

\section{Basic set-ups and holomorphic extendibility of a  reflection function}
 
\bigbreak
In this section, we  first set up some notation which will be used throughout
 the paper. Then we study the holomorphic extendibility of a reflection
function. We will use some notations, which the author adapts from [DP].

To start with, we
let $M_1\subset {\bf C^{n+1}}$, $M_2\subset
{\bf C^{n+1}}$
 be  real analytic hypersurfaces of finite D-type, and
let $f$ be a {\it non-constant} CR mapping from $M_1$ to $M_2$. Assume that
$f$ extends as a holomorphic correspondence to a neighborhood of $p(\in M_1)$.
Our main objective from $\S 2$ to $\S 5$ is to understand under
what circumstances $f$ admits a holomorphic extension near $p$.
Since this is purely a local problem, 
 after  holomorphic changes of coordinates and shrinking the size
 of the hypersurfaces, we can assume  that
$p=f(p)=0$; $M_1$ , $M_2$ are connected and defined, respectively,
 by the functions
of the following form (see [BJT], for example):

$$\rho_1(z,\overline z)=z_{n+1}+\overline {z_{n+1}}+\sum_{j=0}^{\infty}\phi_j(z',\overline {z'})
\left({z_{n+1}-\overline {\zn}\OO 2i}\right)^j,\leqno (2.1)$$
$$\rho_2(z,\overline z)=z_{n+1}+\overline {z_{n+1}}+\sum_{j=0}^{\infty}\psi_j(z',\overline {z'})
\left({z_{n+1}-\overline {\zn}\OO 2i}\right)^j,\leqno (2.2)$$
with $\phi_j(z',0)=\phi_j(0,\overline {z'})=\psi_j(0,\overline {z'})=\psi_j(z',0)\equiv 0$
 and $\phi_j$,$\psi_j$ holomorphic
in $z',\overline {z'}$ over certain fixed open neighborhood of $(0,0)$.
Here we use $z=(z',\zn)$ with $z'=(z_1,\cdots,z_n)$ for the coordinates
in $\bf C^{n+1}$.

Applying a CR extension result in [BT] or [Tr] and making
$M_1$ sufficiently small, we can 
further assume that
there is a domain $D$ which contains $M_1$ as part of its smooth boundary
such that every CR function defined over $M_1$ can be extended holomorphically
to $D$. Hence, without loss of generality, we can assume, in all that follows,
that $f$   is holomorphic in $D$ and continuous up to $D\cup M_1$.

Applying the implicit function theorem to (2.1), (2.2) and shrinking the
size of the hypersurfaces further (if necessary), we see that $M_1$ and
$M_2$ can also  be presented, respectively, by functions
of the following form:
$$\wt{\rho_1 }(z,\overline {z})=\zn+\overline {\zn}+\sum_{j=0}^{\infty}\wt{\phi_j}(z',\overline {z'})
\zn ^j=0,\leqno (2.1)'$$
$$\wt{\rho_2 }(w,\overline {w})=z_{n+1}+\overline {z_{n+1}}+\sum_{j=0}^{\infty}\wt{\psi_j}
(z',\overline {z'})
z_{n+1}^j=0,\leqno(2.2)'$$
where $\wt{\phi_j}$ are determined by $\phi_j$'s,
 $\wt{\phi_j}(z',0)=\wt{\phi_j}(0,\overline  z')=\wt{\psi_j}(z',0)
=\wt{\psi_j}(0,\overline {z'})=0$; and 
$$\wt{\rho_j}=\rho_j h_j.\leqno(2.3)$$
Here $h_j$ is a real analytic function near $M_1$ for $j=1,2$ and
$h_j(0,0)=1$. Complexifying
(2.3), we then 
obtain
$$\wt{\rho_j}(z,\overline  w)=\rho_j(z, \overline  w) h_j(z, \overline w).\leqno(2.3)'$$
Therefore, for $z, w$ close to $0$, 
${\rho_j}(z,\overline  w)=0$ if and only if 
$\wt{\rho_j}(z, \overline  w)=0$.

We now recall the definition of Segre varieties and some related
  notions, which
were first used in [We] (see also [DW] and [DF1])
 for the study of the mapping problems.
 Since $M_2$ is completely symmetric to $M_1$,
we will only introduce notations for $M_1$
 and add $'$ for those corresponding 
to $M_2$. 

We first choose a small neighborhood ${\cal U}_3$ of $0$ such that $\rho_1(z,w)$,
$\wt{\rho_1}(z,w)$, $h_1(z,w)$ are holomorphic  over
${\cal U}_3\times \h{conj}({\cal U}_3)$, 
 and $h_1(z,w)\not =0$ in
${\cal U}_3\times \h{conj}({\cal U}_3)$.
Here and in what follows,
for a subset $A$ in the complex spaces we write
$\h{conj}(A)=\{\overline {z}: z\in A\}$. For each point $w\in {\cal U}_3$,
the Segre variety associated to  $M_1$ and $w$ is defined by
$Q_{w}=\{z\in {\cal U}_3: \rho_1(z,\overline {w})=0\}$. By making ${\cal U}_3$
small, it can be seen that $Q_{w}$ is actually
a complex submanifold  of ${\cal U}_3$.
It is easy to show that there are two small neighborhoods $
{\cal U}_2,\ {\cal U}_1 $  of $0$ such that
the following holds (see [DF1] for more details on this matter): 

\noindent (I.a)
${\cal U}_3\supset
\supset
{\cal U}_2\supset\supset {\cal U}_1\ni 0$.

\noindent (I.b) 
 ${\cal U}_1\sm M_1$ has 
two connected components,
 which are in different sides of $M_1$. Moreover,
they are homeomorphic to the ball. 

\noindent (I.c) For each $w\in {\cal U}_1$,
  ${\cal U}_2\cap Q_w$ is connected.

\noindent (I.d) For a sequence $\{w_j\}\subset {\cal U}_1$ with
$w_j\ra w\in {\cal U}_1$, the limit set of $\{Q_{w_j}\}$, denoted
by $\lim Q_{w_j}$, is $Q_{w}$.

 In what follows, we call $\{{\cal U}_3, {\cal U}_2, {\cal U}_1\}$
a Segre neighborhood system around $0$. It is clear that for any open set 
containing
 $0$, we can always construct inside this set
 a Segre neighborhood system around $0$. Meanwhile, when ${\cal U}_3$
(${\cal U}_2$, respectively)
is fixed, we can arrange ${\cal U}_2$  (${\cal U}_1$,
respectively) as small as we need.

For convenience of the reader, we indicate here 
how (I.c) and (I.d) can be achieved (much more details for related
discussion can be found in work of [DW] and [DF1]).

In the coordinates system where $M_1$ is defined by (2.1), we can choose
a small polydisc $P=P(\delta_1,\cdots,\delta_{n+1})=
\{z:\ |z_j|<\delta_j
(\ j=1,\cdots, n+1)\}(\subset\subset
{\cal U}_3)$. Then, applying the implicit function theorem,
we can see that for each $w$ with $|w|$ small,
 $Q_{w}\cap P$  is contained
in the graph of  certain function , say $z_{n+1}=g(z',\overline w)$, over
the polydisc
$P'=P'(\delta_1,\cdots,\delta_n)=
\{z':\ |z_j|<\delta_j,\ j=1,\cdots, n\}$. Here $g(z',\overline w)$ is holomorphic
for $z, w\in {\cal U}_3$ and $g(z',0)=0$. Now, if we make
 $w$ close enough to the origin,
 then $Q_{w}\cap P$ will be exactly the graph of $g(z',\overline w)$
over $P'$. Since it is biholomorphic to $P'$, it is of course connected.
 Hence, we can
simply take ${\cal U}_2$ to be the polydisc  $P$ chosen in such a manner.
Once ${\cal U}_2$ is chosen, it is then clear that
 we can simply take ${\cal U}_1$
to be a sufficiently small ball around $0$.
 Next, if $w_j\ra w \in {\cal U}_1$,
then we clearly see that $Q_{w_j}\ra Q_{w}$; for
 $g(z',\overline {w_j})\ra g(z',\overline w)$.
This gives (I.d).

We next set up the notation for a reflection mapping ${\cal R}$.

By examining the proof of Lemma 1.1 of [BJT], it is clear that
for each point $p \in M_1$ close to $0$, 
there exists a biholomorphic mapping
$\Phi(\cdot;p)$,
which depends real analytically on $p$,
 such that the following
holds:

\noindent (i) $\Phi(\cdot;0)=\h{id}$ and $\Phi(p;p)=0$.

\noindent (ii)   $\Phi(\cdot;p)$ 
(and thus its inverse) is well-defined over
some fixed open neighborhood of $0$; 

\noindent (iii) 
$\Phi(M_1;p )$  is defined near $0$ by the following
equation:
$$\rho_1(z,\overline z;p)=z_{n+1}+\overline {z_{n+1}}+\sum_{j=0}^{\infty}\phi_j(z',\overline {z'};p)
\left({z_{n+1}-\overline {\zn}\OO 2i}\right)^j,$$
with $\phi_j(0,\overline {z'};p)=\phi_j(z',0;p)\equiv 0$ and $\phi_j$ holomorphic
in $z',\overline {z'}$ over certain fixed open neighborhood of $(0,0)$.

For a small positive $\epsilon_0$, define 
${\cal \Phi}$: $(-\epsilon_0,\epsilon_0)\times M_1\rightarrow {\bf C^{n+1}}$
 by ${\cal \Phi}(t,p)=\Phi^{-1}((0',t);p)$.
Using the facts that ${\cal \Phi}(0,p)=p$ and ${\cal \Phi}(t,0)=(0',t)$,
it can be seen that ${\cal \Phi}$ is an analytic diffeomorphism near
$(0,0)$.

We now  define the reflection mapping ${\cal R}$ as follows:
For $z$ close to $0$, 
if ${\cal \Phi}^{-1}(z)=(t,p)$, then ${\cal R}(z)={\cal \Phi}(-t,p)$.
Clearly, ${\cal R}$ is a diffeomorphism near the origin. By
 the way it was constructed,
one can verify that ${\cal R}|_{M_1\cap O(0)}=\h{id}$, ${\cal R}^2(z)=z$
for $z$ close to $0$,
${\cal R}(D\cap O(0))\subset D^c$,  ${\cal R}(D^c\cap O(0))\subset D$. Here
$D^c={ U}\sm ( D\cup M_1)$ with $U$ a  neighborhood of $M_1$; and
in all that follows we write $O(a)$ for a {\it small}
 neighborhood of $a$, whose
size may be different in different contexts. 

We observe that for a real number $t$ with $|t|$ small, the Segre variety
of $\Phi(M_1;p)$ at $(0,t)$ is an open subset of the 
affine space $\{(z',-t)\}$
near $(0,-t)$. Using the invariant property of Segre varieties under 
biholomorphic mappings ([We] [DW]) and by arranging the size of
$\{{\cal U}_3,{\cal U}_2,{\cal U}_1\}$ suitably, 
one can see that  ${\cal R}(z)\in Q_{z}$
for each $z\in {\cal U}_2$.

Since we assumed that $f$ extends as a correspondence to a neighborhood of $0$,
after shrinking $M_1$ if necessary, we can assume that for each component
$f_j$ of $f$, there is an {\it irreducible} Weierstrass polynomial $${\cal P}_j(z;X)
=X^{N_j}+\sum_{0\le l<N_j}a_{j,l}(z)X^l\leqno (2.4)$$
with $a_{j,l}$ holomorphic over a neighborhood, say $U$, of $M_1$
such that ${\cal P}_j(z;f_j(z))=0$ for $z\in M_1$. Since $f$ is known
 holomorphic on $D$, ${\cal P}_j(z;f_j(z))=0$ for
$z\in D\cup M_1$.
In all that follows, we always make (a) $U\supset\supset {\cal U}_3\cup D$;
(b) $D\supset\supset {\cal U}_3\cap D$; (c) ${\cal R}$ is a diffeomorphism
from ${\cal U}_3$; and (d) ${\cal R}(D\cap {\cal U}_1)\subset D^c\cap
{\cal U}_2$, ${\cal R}(D^c\cap {\cal U}_1)\subset D\cap
{\cal U}_2$. 

As usual, we define the branching locus ${\cal E}$ of $f$ to be
the union of the branching locus of $f_j$ ($j=1,\cdots, n+1)$. That is,
a point $z\in { U}$ is in ${\cal E}$ if and only if for some index
$j$ and
number $X$, one has
 ${\cal P}_j(z;X)=
{\PP {\cal P}_j\OO \PP X}(z;X)=0$. 
Write $E_{M_1}={\cal E}\cap M_1$. Then by the finite type assumption of $M_1$,
one can see that $E_{M_1}$ is a real analytic subset of $M_1$ with codimension
$\ge 2$. Obviously,
 $f$ admits a holomorphic extension across $M_1\sm E_{M_1}$.

We now introduce the following $\lambda$-function associated to the map $f$
and $M_2$ (see [Le], [Pi1], [BJT], and, in particular, [BBR],  for
closely related notions): 
$$G(f(z),\l)=-f_{n+1}(z)-\sum_{j=0}^{\infty}\wt{\psi_j}(f^*(z),\l)
f^{j}_{n+1}(z).
\ \ \ z\in M_1\cup D\leqno (2.5)$$
Here
$\wt{\psi_j}$ are the same as in  $(2.2)'$ and
we write $f^*=(f_1,\cdots,f_n)$.
 Then it is clear that $G(f(z),\l)$ is holomorphic over $D\times { O}_{\l}(0)$
and continuous on
 $(D\cup M_1)\times { O}_{\l}(0)$.  Here and in all that
follows, we put a  subscript to $O(a)$
to emphasize the coordinates used for the space where
$O(a)$ stays. For example, $O_b(a)$ denotes 
a small open neighborhood of $a$
 in the complex b-space, whose size may be different in different contexts.

Our first  lemma in this paper is to extend this function across $M_1$ for
 small $\l$. 
\bigbreak
{\bf Lemma 2.1}: Under the above notations and assumptions, 
$G(f(z),\l)$ extends
holomorphically to ${ O}_z\times { O}_{\l}$. 
\bigbreak
The proof we present here  is to take the differentiation along the boundary,
as did  in [Le], [Pi], 
[BJT], 
[BBR], [BR1],  [Hu], and [BHR1] (in particular, see
[BJT], [BBR], and $\S 2.4.2$ of [Hu] or Proposition
7.1 of [BHR1]).
 However, there is an essential
  difference here.
That is,  our map is not assumed to be smooth.
So, we can only do it almost everywhere.
 To reach the bad points, we jump into the domain and 
  use the hypothesis to control the rate of blowing-up.

 We first observe that when $D\cup M_1$ is pseudo-concave along
$M_1$, then clearly $f$ extends holomorphically across $M_1$. So, 
without loss of generality, we will always
assume the existence of a non empty pseudoconvex piece in $M_1$. 
Furthermore, using the finite type assumption for $M_1$, we can then 
see the existence of strongly pseudoconvex points in $M_1$ 
(see, for example, [BHR1]).
 Since
it is now clear that  strongly pseudoconvex points form a non empty
 open
 subset of $M_1$ and since $f$ extends holomorphically across
$M_1\sm E_{M_1}$, by our non constant assumption of $f$ and
 by a result in [BR3], 
the Jacobian $J_f$ of $f$ is well-defined and non zero on an open
dense subset 
in $M_1$.  

\bigbreak
{\it Proof of Lemma 2.1}:  By the definition of $G(f(z),\l)$
and using the assumption that $f(M_1)\subset M_2$,
 we have
$$\overline {f_{n+1}(z)}=G(f(z),\overline {f^*(z)})\ \ \ \h{for}\  z\in M_1.\leqno (2.6)$$

Shrinking the size of $M_1$ if necessary, we
can choose a basis  $\{{\cal L}\}_{j=1}^{n}$ of 
 $\h{T}^{(1,0)}M_1$, whose coefficients are real analytic in $z$.
 As for the smooth mappings case (see [BJT], [BR1]),
we apply $\overline {\cal L}_j$ to (2.6),
to obtain
$$\overline {\cal L}_j{\overline  f_{n+1}}(z)=\sum_{l=1}^{n}{\PP G\OO\PP \l _{l}}\overline {\cal L}_j
\overline {f}_l,\ \ \ \h{for} \ z\not\in E_{M_1}.$$
Let $J=\h{det}(\overline {\cal L}_j\overline {f_l})_{1\le j,l\le n}$ and let
$\cal J$ be the matrix 
$(\overline {\cal L}_j\overline {f_l})_{1\le j,l\le n} $
, which are well-defined 
 on $M_1\sm E_{M_1}$.
 We claim that $J\not\equiv 0$ on  $M_1\sm E_{M_1}$.
 To see this,   by the above observation,
we can assume that
there is 
a strongly pseudoconvex point $p_0\in M_1$ such that $f$ extends
biholomorphically 
to a neighborhood of
 $p_0$ and maps $p_0$ to a strongly pseudoconvex point $q_0$ of $M_2$. 
 Now, seeking a contradiction suppose $J(z)\equiv 0$
near $p_0$. Let ${\cal J}_j(z)=({\cal L}_jf_1,\cdots, {\cal L}_jf_n)$.  
By passing to a nearby point,
we can assume, without loss of generality, that
the semi-continuous function $d(z)$, which is defined as
the complex dimension of  the vector space $V(z)$   spanned by
$\{{\cal J}_j(z)\}$, takes a local maximal value, say $k$ $(<n)$, at $p_0$.
Here
 z is close to $p_0$ but stays in $M_1$. Moreover,
we can assume that
 $\{{\cal J}_1(p_0),\cdots,{\cal J}_k(p_0)\}$ forms a basis of $V(p_0)$.
Now, by 
the local maximality of $d(z)$ at $p_0$ and by
the real analyticity of ${\cal J}_j(z)'s$ with respect to $z$,
 it is clear, from some simple arguments in linear algebra,
 that near $p_0$ there exist real analytic functions
$a_j(z)$ ($j=1,\cdots,k$)
 with $z\in M_1$ sufficiently close to $p_0$ such that
${\cal J}_{k+1}(z)=\sum_{j=1}^{k}a_j(z){\cal J}_j(z)$
 for $z (\in M_1)\approx p_0$. Write $T(z)={\cal L}_{k+1}(z)-
\sum_{j=1}^{k}a_j(z){\cal L}_j(z)$. Then $T(z)$ annihilates $f^*(z)$
for $z (\in M_1)\approx p_0$. Since $T(p_0)\not =0$ and
$M_1$ is strongly pseudoconvex at $p_0$, it follows that
$(d\rho_1, [T,\overline {T}])(p_0)\not = 0$,
 and thus $\{{\cal L}_1,
\cdots, {\cal L}_n,[T,\overline {T}]\}$ 
forms a basis of $T_{p_0}M_1$.
 Write $N(z)$ for the $(1,0)$-component of $[T,\overline {T}](z)$.  
Noticing
that $f$ is holomorphic near $p_0$ and
 $[T,\overline  T]f^*=T((\overline {T}f^*))-\overline {T}(Tf^*)=0$,
 it is clear that $N(p_0)$ annihilates $f^*$, too. Observe that
$\{{\cal L}_1,\cdots, {\cal L}_n, N\}(p_0)$ forms a basis of
$T^{(1,0)}_{p_0}{\bf C^{n+1}}$. Therefore, making use of the fact $N(p_0)(f^*)=0$ and by a simple linear algebra argument, it follows that
 $J_f(p_0)$ is given by some constant times $\overline {J(p_0)}$ and thus equals to $0$.
 This
contradicts our choice of $p_0$.

Now, we let  
  $Z^*=\{z\in M_1\sm E_{M_1}: \  J(z)=0\}$.
Then, since $M_1\sm E_{M_1}$ has to be connected,
we see that  $M_1\sm Z^*$ is a dense open subset in  $M_1$.
In fact, with a little  more effort, we will next show that
$Z^*$ is contained in a real analytic subset of codimension
at least 1. 
Indeed, notice that $J(z)=J(z,\overline z,\overline {Df})$ can be written as a
polynomial in $Df=({f_1}'_{z_1},\cdots,{f_n}'_{z_{n+1}})$ with
coefficients holomorphic in $(z,\overline z)$ for $z$ in a small neighborhood
of $M_1$. Since $(2.4)$ annihilates $f_j$, one can easily show that
 each component of $Df$ can be presented as a rational
 function of $f$ with  coefficients  holomorphic near $M_1$.
Thus, by  some standard algebra arguments involving
 symmetric functions, one
can conclude that for some irreducible polynomial in $X$:
${\cal P}(z,w;X)=\sum_{j=0}^{N^*}A_j(z,w)X^j$, it holds that
${\cal P}(z,\overline z;J(z,\overline z, \overline {Df(z)}))=0$ for $z \in M_1\sm E_{M_1}$. Here
$A_j(z,w)$ are holomorphic for $z$ and $\overline  w$ in a small neighborhood
of $M_1$. Since $A_j(z,\overline z)=0$ on $M_1$ if and only if $A_j(z,\overline z)=\rho_1(z,
\overline z)h_j(z,\overline z)$ for some real analytic function
 $h_j$ near $M_1$, after doing some cancellation, 
 we may assume that not all $A_j(z,\overline  z)$
vanish on $M_1$.
 Now, since $J(z)\not = 0$ on a dense open subset,
the zero set of $J(z)$ is obviously contained in the
zero set of the coefficient $A_{j^*}(z,\overline z)$, where  $A_{j^*}(z,\overline z)\not 
\equiv 0$ but $A_{j}(z,\overline  z)\equiv 0$  for $z\in M_1$ and for $j<j^*$.
 Hence $Z^*$ is contained in the proper real analytic
 subset of $M_1$
consisting of $E_{M_1}$ and
 the zeros of $A_{j^*}(z,\overline z)$ in $M_1$.

 Away from $Z^*$, we now have
$$\left({\PP G\OO\PP \l _1},\cdots, {\PP G\OO\PP \l _n}\right)^t
={\cal J}^{-1}(z,\overline {z},\overline {Df})(\overline {\cal L}_1\overline {f}_{n+1},\cdots, \overline {\cal L}_n\overline {f}_{n+1})^t.
\leqno (2.7)$$
Writing $(2.7)$ as $n$ scalar equations,  
applying $\overline {\cal L}_j$ to each of them, and keeping doing in this manner,
 we see, by induction,
 that for each multi-index ${\bf \a}=(\a_1,\cdots,\a_n)$; there are
 two holomorphic functions $g_{\bf \a}^{(1)}$ and $g_{\bf \a}^{(2)}$ in
 $$(z,\overline {z},\overline {f},\cdots,\overline {D^{\beta}f},\cdots)$$ ( where
 $\beta\le |\a|$)
 such that for each $z\in M_1\sm Z^*$, one has
 $g_{\a}^{(2)}(z,\overline z,\cdots,\overline {D^{\beta}f},\cdots)\not =0$ and
$$D^{\a}_{\l}G(f,\overline {f^*(z)})={  
g_{\a}^{(1)}(z,\overline z,\cdots,\overline {D^{\beta}f},\cdots) \OO
g_{\a}^{(2)}(z,\overline z,\cdots,\overline {D^{\beta}f},\cdots)}.$$
Here $D^{\b}$ denotes the vector formed by
all derivatives of  $f$ with order $\b$.
 We remark that $g^{(j)}_{\a}$ is actually a polynomial in
$(\overline {Df},\cdots,\overline {D^{\b}f},\cdots)$ ($j=1,2$) with coefficients real analytic in $z$ for $z$ close to $M_1$.
By passing to the limit, we see that
the function ${g^{(1)}_{\a}\OO g^{(2)}_{\a}}$ has a continuous extension
to $M_1$, which we will denote by $h_{\a}(z,\overline  z,\cdots,\overline { D^{\b}f},\cdots, )$.
 Notice also 
that for $\|w\|, \|\l\|\ll 1$, there exits a large constant
$R$ so that
$|D^{\a}_{\l}G(w,\l)|\ale \a !R^{|\a|}$. 


Denote all the solutions of (2.4) by $\eta_{j,l}(z)$ $(j=1,\cdots,n+1,\
l=1,\cdots,N_j)$. 
We now can easily find a totally real submanifold $M'\subset M_1$ of dimension
$n+1$ such that $0\in M'$
and $M'\cap Z^*$ is contained in a real analytic
 subset of  $M'$ with codimension at least 1 (see, for example,
 [Hu] or [BHR1]).
 Choose ${\cal H}$:
$({\bf C^{n+1}},0)\ra ({\bf C^{n+1}},0)$, a germ of biholomorphism which
sends a small open neighborhood of $0$ in  ${\bf R}^{n+1}$
into $M'$ and maps the standard wedge $W^+$ into $D$,
$W^-$ into $D^c$, respectively. Here the notation $D^c$ is as before. 

Consider the following equation in $X$:
$$\prod_{l_1,\cdots,l_{2n+1}}
 \left(X-{1\OO \a!}(D^{\a}_{\l}G(\eta_{1,l_1}\circ
{\cal H}(z),\cdots,\eta_{n+1,l_n}\circ
{\cal H}(z), \cdots, \overline {\eta_{n,l_{2n+1}}\circ {\cal H}(\overline {z})})\right)$$
$$=X^N+\sum_{j<N}c_j(z)X^j=0,$$
where $l_j$ and $l_{n+1+j}$ run from $1$ to $N_j$ for $j\le n+1$, and
$N=N_1^2\cdots N_n^2N_{n+1}$. 
$c_j(z)'$s can be seen to be the symmetric functions of
some Weierstrass polynomial equations and hence can be seen to be
holomorphic near $0$. Moreover,
one can obtain the Cauchy estimates $|c_j(z)|\ale R^{N|\a|}$.
Now, we notice that
$$\wt{h}_{\a}(z)={1\OO \a!}h_{\a}({\cal H}(z),\overline {{\cal H}(\overline  z)},\cdots,
\overline {D^{\beta}f({\cal H}(\overline {z}))},\cdots)$$
 is a solution of
the above equation for almost all $z\in {\bf R}^{n+1}$ near $0$. Meanwhile,
it clearly extends to a meromorphic function to $W^-\cap O_z(0)$;
for $f$ is holomorphic over $D$.
 Thus, by the uniqueness
of holomorphic functions [Pi2], it follows that
$$(\wt{h}_{\a}(z))^N+\sum_{j<N}c_j(z)(\wt{h}_{\a}(z))^j=0,$$
for all $$z\in W^{-}\cap { O}_z(0)\sm \{\h{the singular set of }\ 
\wt{h_{\a}}\ \h{in}\ W^-\cap O_z(0)\}.$$
In particular,
we see that
$\wt{h}_{\a}(z)$ is bounded.  Using
 the Riemann extension theorem, we conclude that
$\wt{h_{\a}}$
 extends holomorphically to $ W^{-}\cap O_z(0)$.
Moreover, we   have
$|\wt{h}_{\a}(z)|
\ale R^{N|\a|}$
for $z\in W^{-}\cap { O}_z(0)$
; for its coefficients have the same sort of estimates.

Next,  fix a small open subset $U$ containing $0$. 
Let  $$\phi^-_{\a}(z)={1\OO \a !}\left (D^{\a}_{\l}\sum_{\b} 
 \wt{h}_{\b}(z)\left(\l-\overline {f^*\circ{\cal H}(\overline  z)}
\right)^{\b}\right )_{\l=0}$$ for
  $z\in W^-\cap U$ 
and let $\phi^+_{\a}(z)={1\OO \a !}D^{\a}_{\l}G(f\circ {\cal H}(z),\l)|_{\l=0}$
for  $z\in W^+\cap U$. Then it can be seen that $|\phi^+_{\a}(z)|,\
 |\phi^-_{\a}(z)|
\ale {R'}^{|\a|}$ for some large constant $R'$, where $z$ stays in
 their defining
regions, respectively. Notice that
  $\phi_{\a}^+$ matches up with $\phi_{\a}^-$
 over $U\cap
{\bf R}^{n+1}\sm \{\h{ a thin real analytic subset}\}$. 
The classic edge of the wedge theorem then indicates
that $\phi^+_{\a}$ extends to a holomorphic function $\phi_{\a}$ defined over
 some sufficiently small neighborhood $U'$ of $0$,
 whose size depends only on the size of the wedges where
$\phi^{\pm}_{\a}$ are defined, and therefore is
 independent of $\a$ (see [Vad] or [BHR1], for example).
 Moreover, $U'$ can be filled in
by analytic disks  with boundary staying in  the closure of
 $(W^{+}\cup W^{-})$ ([Val]). So, the maximal principle tells that
$\phi_{\a}$ has the same kind of Cauchy estimates as $\phi_{\a}^+$ and $\phi_{\a}^-$ do.
Now, $\left(\sum_{\a}\phi_{\a}\l ^{\a}\right)\circ {\cal H}^{-1}$
 clearly gives the holomorphic
extension for $G(f(z),\l)$ to ${ O}_z(0)\times { O}_{\l}(0)$.
 The proof
of Lemma 2.1 is now complete. $\endpf$
\bigbreak

{\bf Remark 2.2}:
(a) As an immediate application of Lemma 2.1,
by letting $\l=0$ in (2.4), we conclude that $f_{n+1}(z)$ admits a holomorphic
extension near $0$. Thus,  ${\cal P}_{n+1}(z;X)=X-f_{n+1}(z)$.


(b) From the proof of Lemma 2.1, one can see that the same result holds
 whenever $f$ extends to one side and
  $M_1\sm {\cal E}$ contains a strongly pseudoconvex point where
$J_f$ is not 0. In particular, when $f$ extends to one side,
 this can be applied to the case
in which $M_1$ and $M_2$ are
not Levi-flat in $\bf C^2$.

(d) In the following sections, we always assume that $G(f(z),\l)$
is holomorphic over the closure of ${\cal U}_3\times { O}_{\l}(0)$
 
  
\section{Branches  of $\cal F$ and Segre varieties}

In this section, we start to study the branching property of the polynomial
equations defining $f$. The argument will be based on 
the extensive use of Segre
varieties and Lemma 2.1. 
Before proceeding, we set up the following notations:

Let $\cal E$ still be as before.
For each $z\in {\cal U}_3$, we write
${\cal F}(z)$ for the multiple-valued map formed by extending
  $f$ from
${\cal U}_3\cap D$. More precisely, when $z\not\in {\cal E}$,
 $\wt{f}(z)\in {\cal F}(z)$ if there
exists a path $\gamma\subset {\cal U}_3\setminus {\cal E}$
 with $\gamma (0)\in D$,
 $\g(1)=z$, and a continuous section ${\cal G}$  from $[0,1]$
 to the sheaf of germs of holomorphic mappings to
$\bf C^{n+1}$, denoted by $({\cal O}^{n+1},\pi, {\bf C^{n+1}})$,
 such that
${\cal G}(0)$ coincides with the germ of $f$ at $\g(0)$
and $g(\g(1))=\wt{f}(z)$,  where $g$ is a representation of ${\cal G}(\g(1))$,
i.e, the germ of $g$ at $\g(1)$ is
${\cal G}(\g(1))$.  When $z\in {\cal U}_3\cap {\cal E}$, 
 we let ${\cal F}(z)$ be the
 set of the limit
points of all possible sequences $\{{\cal F}(z_j)\}$ with $z_j\not \in 
{\cal E}$ and
$z_j\rightarrow z$.
For each $z\in {\cal U}_2$, we write
${\cal F}(Q_z)=\cup_{w\in Q_z\cap {\cal U}_2}{\cal F}(w).$
  We also notice that for each point in ${\cal F}(z)$, its $m-$th component
is a solution of the equation ${\cal P}_m(z;X)=0$ ((2.4)). From this, it follows
that when $z$ is close to $0$, then ${\cal F}(z)$ is close to the origin, too.
 Moreover, 
Lemma 2.1 now indicates
that $G({\wt{f}}(z),\l)$ is independent of the choices 
of $\wt{f}(z)\in {\cal F}(z)$. So, in what follows, we 
use the notation $G({\cal F},\l)$ to replace $G(f(z),\l)$.

For $w\in {\cal U}_2$, we denote by $A_w$ the set $\{\o\in {\cal U}_2: Q_w\cap
{\cal U}_2=Q_{\o}\cap {\cal U}_2\}$. We  mention that
$A_w$ is invariantly attached to the hypersurafce $M_1$, and will
play an important role in our proof of Theorem 1.1. 
 By  our finite D-type assumption
for $M_1$ and 
 by shrinking ${\cal U}_3$, we can make $A_0=\{0\}$
and $A_z$ finite for $z\in {\cal U}_2$
(see [DW], [BJT], and [DF1]). Similarly, by making ${\cal U}'_3$ small,
we can assume that $A_0'=\{0\}$ and $\#A'_{w}<\infty$ for
$w\in{\cal U}'_2$. For more properties of $A_w$, we refer the reader to
([DF1], Lemma 1).

Before proceeding to the proof of  the first  lemma in this section,
 we
mention that  by shrinking ${\cal U}_2$ if necessary,
we can assume that ${\cal F}({\cal U}_2)\subset {\cal U}_1'$.
Meanwhile, it is clear that ${\cal F}(0)=0$.
 We also retain notations which have been 
set up so far.
\bigbreak
{\bf Lemma 3.0}
 (i)  $\# {\cal F}(z)$ is constant for
$z\in {\cal U}_3\sm {\cal E}$.

\noindent  (ii) The limit set of ${\cal F}(z_j)$ is 
${\cal F}(z_0)$ if $z_j\ra z_0 (\in {\cal U}_3)$.

\noindent
 (iii) In case $n=1$, 
${\cal F}(z)=\{(X,f_2(z)):\ {\cal P}_1(z;X)=0\}$.

\noindent (iv) For $w_1, w_2\in {\cal U}_1$, if the variety
$Q_{w_1}\cap Q_{w_2}\cap {\cal U}_2$ contains a point of complex dimension
$n$, then 
$Q_{w_1}\cap {\cal U}_2= Q_{w_2}\cap {\cal U}_2$.  Thus
 $A_{w_1}=A_{w_2}$.

\bigbreak
{\it Proof of Lemma 3.0}:
 (i) is an immediate
consequence of what is called the monodromy theorem. (iv) is a consequence
of (I.c) and the fact that 
two connected closed complex submanifolds of a domain
coincide if they have an open subset in common.
(iii) follows easily from the irreducibility assumption of ${\cal P}_1$
and the following standard argument:
For $z\not \in {\cal E}$, write ${\cal F}(z)=\{(g_j(z),
 f_2(z))\}_{j=1}^{k}$. Then (i) indicates that $k$ is constant. Now,
the Riemann extension theorem tells that $\wt{{\cal P}_1}=\prod_{j=1}^{k}
(X-g_j(z))$ is a Weiestrass polynomial in $X$. It, of course,
 divides ${\cal P}_1$ by the previous observation.
So, from the irreducibility of ${\cal P}_1$, it follows that
${\cal P}_1=\wt {\cal P}_1$, and this clearly completes the proof of (iii).
 So, it now suffices for us to give the proof of (ii) to complete
the proof of Lemma 3.0.

When $z_0\not \in {\cal E}$, then ${\cal F}(z)$ can be stratified into $k$
holomorphic branches near $z_0$ with $k=$the generic counting number
of ${\cal F}(z)$. That is, we can write 
${\cal F}(z)=\{g_1(z),\cdots, g_{k}(z)\}$
with $g_j$'s holomorphic near $z_0$. Now, when $z_j$ is close enough to
$z_0$, then ${\cal F}(z_j)=\{g_1(z_j),\cdots,g_k(z_j)\}$. From this,
it is easy to see that $\lim {\cal F}(z_j)= {\cal F}(z_0)$. When
$z_0\in {\cal E}$, by passing to a nearby point, we can assume, without
loss of generality, that all $z_j$'s are not in ${\cal E}$. Now, for
another sequence $\{w_j\}(\subset {\cal U}_3\sm {\cal E})$ with $w_j\ra z_0$,
we first find  a family of 
curves $\{\g_j\}$ such that $\|\g_j-z_0\|\ra 0$ and
$\g_j(0)=z_j$, $\g_j(1)=w_j$. Still write ${\cal F}(z)=\{g_1(z),\cdots,
 g_k(z)\}$ for $z\not\in {\cal E}$. Then for each $l$, there is a continuous
map $I_l(t)$ from $[0,1]$, which maps $t$ to a point in ${\cal F}(\g_j(t))$
 such that
$I_l(0)=g_{l}(z_j)$ and $I_l(1)=g_{\sigma_{z_j} (l)}(w_j)$, where 
$\sigma_{z_j}$ is
a permutation of $\{1,\cdots,k\}$. Moreover, it is clear that
 the m-th component of $I_l(t)$
is a solution of ${\cal P}_m(\g _j(t);X)=0$, as observed above. Now, by Lemma 2.5 of [Mal]
and our choices of $\g_j$, it follows that $\|g_l(z_j)-g_{
\sigma_{z_j}(l)}(w_j)\|\ra 0$.
Therefore, we see that $\lim {\cal F}(z_j)=\lim {\cal F}(w_j)$, which
must be ${\cal F}(z_0)$ by the definition. $\endpf$
\bigbreak

{\bf Lemma 3.1}. 
After shrinking $ {\cal U}_1$, if necessary,
  then for each $z\in {\cal U}_1$,
it holds that
  ${\cal F}(z)\subset A'_{\wt{f(z)}}$ and ${\cal F}(Q_z)\subset Q'_{\wt f(z)}$
 for any 
 $\wt{f}(z)\in {\cal F}(z)$.
\bigbreak
{\it Proof of Lemma 3.1}: We would like to mention that by a
simple unique continuation argument and by using the invariant property
 of Segre varieties, one can easily show that for each
nice branch $\wt{f}$ with $\wt{f}|_{D}=f$,
 it holds that $\wt{f}(Q_z)\subset Q_{\wt {f(z)}}$ (see also [DF2],
 for example). 
The main idea of the proof of the lemma is to use the fact that
$\wt{\rho}_2({\cal F}(z),\wt{f}(\o))$ is single-valued for each fixed $\o$,
 by  the
above established $\l$-function.
 
 We let $z\in {\cal U}_1\cap D$,
 and assume that
$z,\ z^*\not\in {\cal E}$. 
Here, we use $z^*$ to denote
the reflection point of $z$, i.e, $z^*={\cal R}(z)$.
We also choose a simply connected smooth curve $\gamma : [0,1]\ra {\cal U}_3$
 with
$\gamma \subset
{\cal R}({\cal U}_1\sm D)\sm {\cal E}$ and
$\g(1)=z^*$, $\g\cap M_1=\{\g(0)\}$. Here, when there is no confusion arising,
we also use  the  letter $\gamma$ to denote its image set.  Moreover,
we assume that $\g$ intersects $M_1$ transversally at $\g(0)$.
Let $\wt{\g}={\cal R} (\g) \cup \g$. Then $\wt{\g}$ is still a simply 
connected  curve in ${\cal U}_2$ by  the discussions in $\S 2$.
 Thickening $\wt{\g}$ suitably,
 we can then obtain a simply connected
domain, which we will denote by $O (\subset {\cal U}_2)$.
 Now, we  can well define a holomorphic
map $\H{f}$ from $O$, which coincides with $f$ on $O\cap D$. Let
$M_{1c}=\{(z,\o)\in {\cal U}_3 \times \h{conj}({\cal U}_3): 
\wt{\rho_1}(z,\o)=0\}$ and
consider $M^*_{1c}=M_{1c}\cap
 \{O\times\h{conj}(O)\}$. 
 Now, the following
function is well-defined and holomorphic 
over $M^*_{1c}$:
$$\Xi(z,\o)=\overline {f_{n+1}(\overline  \o)}-G({\cal F}(z),\overline {\H{f}(\overline \o)})
\left (=\wt{\rho_2}({\cal F}(z),\overline {\H f (\overline {\o})})\right )$$
which vanishes on $M_1^*=\{(a,\overline  a):a\in  M_1\cap {\cal U}_3\}$. Notice that
$M_1^*$ is  a totally real subset in $M^*_{1c}$ with maximal dimension,
we conclude from the uniqueness property of holomorphic functions [Pi2] that
$\Xi(z,\o)\equiv 0$ in the union, denoted by
$M_{1c}^{**}$,
 of the connected components
 of $M_{1c}^*$ which have non-empty intersections
with
$M_1^*$. By our choice of $\wt{\g}$,
 we see that $(z,\overline {z^*}),\ (z^*,\overline {z})\in M^{**}_{1c}$; for
$(\g(t),\overline {{\cal R}(\g(t)}))\in M^*_{1c}$ 
$({{\cal R}(\g(t)}),\overline {\g(t)})\in M^*_{1c}$ ($t\in [0,1]$)  
connect them
to $M^*_1$, respectively.
Therefore,
$(z,\h{Conj}(Q_{z}\cap O(z^*)))\subset M_{1c}^*$,
$(z^*,\h{Conj}(Q_{z^*}\cap O(z)))\subset M_{1c}^*$ by the definition of
Segre varieties.
 This implies that 
$G({\cal F}(z),\overline {\H{f}(\o)})=\overline {f_{n+1}(\o)}$
for $\o\in Q_{z}\cap O(z^*)$;
and 
$G({\cal F}(z^*),\overline {\H{f}(\o)})=\overline {f_{n+1}(\o)}$ for $\o\in Q_{z^*}\cap O(z)$. 
 Equivalently, we have
$\wt{\rho}_2({\cal F}(z),\overline {\H{f}(\o)})=0$
for $\o\in Q_{z}\cap O(z^*)$;
and 
$\wt{\rho}_2({\cal F}(z^*),\overline {\H{f}(\o)})=0$ for $\o \in Q_{z^*}\cap O(z)$.

From a basic fact of Segre varieties: $w\in Q_z$ if and only if $z\in Q_w$;
 it now follows that  
$\H{f}(Q_{z}\cap O(z^*))\subset \cap Q'_{\wt{f(z)}} $
 with $\wt{f}(z)\in {\cal F}(z)$; and
${\H f}(Q_{z^*}\cap O(z))\subset \cap Q'_{\wt{f}(z^*)}$ 
with ${\wt f}(z^*)\in {\cal F}(z^*)$. 

By slightly perturbing $z$ if necessary, we assume momentarily
 that $J_{\H{f}}\not =
0$ near $z^*$. Now, since $\H{f}(Q_z\cap O(z^*))
\subset Q'_{\wt{f}(z)}$
for any $\wt{f}(z)\in {\cal F}(z)$, and since
 each of them is a connected complex
submanifold of dimension $n$ near $\H{f}(z^*)$,
  all these submanifolds therefore coincide near $\H{f}(z^*)$.
 Hence, it follows from Lemma 3.0 (iv) that
all $Q'_{\wt{f}(z)}\cap {\cal U}_2$
 are the same. So, 
${\cal F}(z)\subset A'_{\wt f(z)}$ for any given $\wt{f}(z)\in {\cal F}(z)$.

Assume also that $J_{f}(z)\not =0$. 
In a similar manner, we then also see that
${\cal F}(z^*)\subset A'_{\wt{f}(z^*)}$ for any given $\wt{f}(z^*)\in {\cal F}(z^*)$.

Next, we note that ${\wt \rho_2}({\cal F}(z),\overline {\H f}(\o))=0$
if and only if $\rho_2({\cal F}(z),\overline {\H f}(z))=0$. Let $q\in
{\cal F}(z)$. Then, by what we just obtained and by the reality of $\rho_2$,
we have $\rho_2(\H{f}(\o),\overline {q})=0$ and therefore,
 $ \wt{\rho_2}(\H{f}(\o),\overline {q})=0$ for $\o\in Q_{z}\cap O(z^*)$.
Now, since Lemma 2.1 indicates that 
 $\wt{\rho_2}({\cal F}(\o),\overline {q})$ is well defined and holomorphic 
for $\o\in {\cal U}_3$, we conclude that
 $\wt{\rho_2}({\cal F}(\o),\overline {q})=0$ for $\o\in {\cal U}_2\cap Q_z$.
This implies that ${\cal F}(Q_z)\subset Q'_q$.

In a similar manner, we can also show that 
 ${\cal F}(Q_{z^*})\subset Q'_{q^*}$ for any $q^*\in {\cal F}(z^*)$.

Now, we will complete the proof of   
 Lemma 3.1 by passing to the limit.

Let $z$ close to $0$  be such that either  $ z\in {\cal E}\cup {\cal R}({\cal E})\cup M_1$ or
 $J_f(z)=0$, or $J_f(z^*)=0$.
 By what we did above, we can find a sequence
$\{z_j\}$ with $z_j\ra z$ so that
${\cal F}(z_j)\subset A'_{\wt{f}(z_j)}$ for any $\wt{f}(z_j)\in
{\cal F}(z_j)$. In the other words,
 $Q'_{\wt{f}(z_j)}=Q'_{\H{f}(z_j)}$
for any $\wt{f}(z_j), \H{f}(z_j)\in {\cal F}(z_j)$. By Lemma
(3.0) (ii) and (I.d), it then follows that
 $Q'_{\wt{f}(z)}=Q'_{\H{f}(z)}$
for any  $\wt{f}(z), \H{f}(z)\in {\cal F}(z)$, i.e,
 ${\cal F}(z)\in Q'_{\wt f(z)}$ for any ${\wt f}(z)\in {\cal F}(z)$.
Similarly, we also have ${\cal F}(Q_z)\subset Q'_{\wt f(z)}$ for
any ${\wt f(z)}\in {\cal F}(z)$.

Hence, after shrinking ${\cal U}_1$ one more time if necessary, we see 
the proof of Lemma 3.1.

 $\endpf$
\bigbreak
{\bf Remark 3.2} (a)
By Lemma 3.1, we can now well-define
 $Q'_{{\cal F}(z)}$ to be $Q'_{q}$
and $A'_{{\cal F}(z)}=A'_q$ for some $q\in {\cal
F}(z)$. Then Lemma 3.1 can be written as ${\cal F}(z)\subset A'_{{\cal F}(z)}$
and
 ${\cal F}(Q_z)\subset Q'_{{\cal F}(z)}$.

(b) As an application of Lemma 3.1,
we conclude  that when $A'_{w}$ is just a single point for $w\approx 0$, then
$f$ extends holomorphically near $0$. This is the case when the target
point is a 
non-degenerate point or has some special bi-type property.
Define ${\cal A}'$ by sending each point $w$ to $A'_w$. One can  see
that if $f$ does not allow holomorphic extension, then ${\cal A}'$ branches
at $0$ (we will make this more precise in $\S 5$).
 It is this fact that links the branching points 
of ${\cal F}$ with the singular
points of  ${\cal A}'$, which will be the key observation
for the proof of Theorem 1.1.

More specifically,
  let $D_1=\{(z_1,z_2): \ |z_1|^2+|z_2|^2<1\}$ and $D_2=\{(z_1,z_2):\
|z_1|^4+|z_2|^2<1\}$. Let ${\cal G}=(\sqrt{z_1},z_2)$ be 
the multiple-valued
map from $D_1$ to $D_2$. Then the branching locus of ${\cal G}$ is
given by $Z=\{(0,z_2)\}$, and 
${\cal G}(Z)=\{(0,z_2)\}$. We observe that ${\cal G}(Z)$ is exactly the
branching locus of the ${\cal A}'$-map of $\PP D_2$.

\bigbreak
{\bf Lemma 3.3}: (i) Let $M_2$ be as given in (2.2).  Let $\O^{\pm}$ 
 be defined by
$\pm{\rho_2}<0$ respectively.  Write $\H{n}^{\pm}=\{b=
(0',b_{n+1}): \pm b_{n+1}<0\}$ with $|b_{n+1}|$ small. 
Then $A'_b\subset \O^{\pm}$ when $b\in \H{n}^{\pm}$, respectively.

\noindent (ii) If for some $q\in {\cal F}(z)$, it holds that $q\in M_2$,
then ${\cal F}(z)\subset M_2$. Hence, ${\cal F}(M_1\cap {\cal U}_1)\subset M_2$.

\noindent (iii) $f_{n+1}(z)=z^k_{n+1}g(z)$ for some holomorphic function
$g(z)$ defined near $0$.

\noindent (iv) $\{{\cal F}^{-1}(0)\}\cap {\cal U}_2=\{0\}$
 after shrinking ${\cal U}_2$, i.e
, $\wt{f}(z)\not =0$ for any $z(\in {\cal U}_2)\not =0$
 and $\wt{f}(z)\in
{\cal F}(z)$. Moreover ${\cal F}^{-1}(w)\cap {\cal U}_2$
 is a finite set for any $w\approx 0$; and for any analytic variety ${ E}$
passing through $0\in {\cal U}_1$, ${\cal F}({ E})$
also gives the germ of an analytic variety  at $0\in {\cal U}'_1$
with $\h{dim}_0 {\ E}=\h{dim}_0({\cal F}({ E}))$. 

\noindent (v) After shrinking ${\cal U}_1$ if necessary,
then for any $z\in {\cal U}_1$, either ${\cal F}(z)\subset \O^+$
or ${\cal F}(z)\subset \O^-$, or ${\cal F}(z)\subset M_2$.
\bigbreak
{\it Proof of Lemma 3.3}: 
(i): Let $\eta=(\eta',\eta_{n+1})\in A'_{b}$, i.e, $Q'_{ \eta}\cap {\cal U}'_2
=Q'_{b}\cap {\cal U}'_2$. 
We see that
$$\{(w',w_{n+1})\in {\cal U}'_2: w_{n+1}+\overline {\eta_{n+1}}+\rho^*_{2}(w',\overline \eta', {w_{n+1}-\overline {
\eta_{n+1}}\OO\ 2i})=0\}$$
$$
=\{(w',w_{n+1})\in {\cal U}'_2: w_{n+1}+\overline {b_{n+1}}+\rho^*_{2}(w',\overline b', {w_{n+1}-\overline {
b_{n+1}}\OO\ 2i})=0\}.
,$$
where $\rho^*=\sum_{j\ge 0}\psi_j(w',\overline {w'})(\h{Im}w_{n+1})^j$.
Letting $w'=0$, we see that $\eta_{n+1}=b_{n+1}$. Meanwhile, 
since
$\{(w',-\overline {\eta_{n+1}})\}\cap {\cal U}'_2= Q'_{\eta}$,
 it follows that
$\rho^*_{2}(w',\overline {\eta'},i\overline {\eta_{n+1}})=0$ for
all $w'$ with $(w',b_{n+1})\in {\cal U}'_2$. 
We see, in particular, that
$\rho^*_{2}(\eta',\overline {\eta'},i\overline {\eta_{n+1}})=0$. Therefore,
${\psi_0}(\eta',\overline {\eta'})=-\sum_{j>0}^{\infty}{\psi_j}(\eta',\overline {\eta'})
(i\overline {\eta_{n+1}})^j=o(|\eta_{n+1}|)$.
Now, $\rho_2(\eta,\overline {\eta})=\eta_{n+1}+\overline {\eta_{n+1}}+
\rho^*_{2}(\eta',\overline {\eta'}, {\eta_{n+1}-\overline {
\eta_{n+1}}\OO\ 2i})$ = $\eta_{n+1}+\overline {\eta_{n+1}}
+{\psi_0}(\eta',\overline {\eta'})+o(|\eta_{n+1}|)
=2\eta_{n+1}+o(\eta_{n+1})$, which is positive when $b_{n+1}=\eta_{n+1}>0$
and small; and negative when $b_{n+1}<0$.
This gives the proof of part (i).

(ii). Let $q\in M_2\cap {\cal F}(z)$. Then $q\in Q'_{q}=Q'_{\wt f(z)}$
for each $\wt{f}(z)\in {\cal F}(z)$. Thus
$\wt{f}(z)\in Q'_{q}=Q'_{\wt f(z)}$,  by Lemma 3.1.
 So, $\wt{f(z)}\in M_2$ (see [DF1]).

(iii). By Lemma 3.1,
 it follows that ${\cal F}(Q_0)\subset Q'_{{\cal F}(0)}=
Q'_{0}$. Thus, we see that $f_{n+1}(z',0)\equiv 0$.
 Since $f_{n+1}$ is holomorphic
near $0$, we conclude that $f_{n+1}(z',z_{n+1})=z_{n+1}^kg(z)$ for some
holomorphic function $g$ near $0$ and  some positive integer $k$.

(iv). The  following argument will be based on some  facts on
 local analytic covers and the elimination theorem (see [Chi], [Wh], or
[GR]).

We start by setting up some notations:

Let ${\wt V}=\{(z, w)\approx (0,0): w_{n+1}=f_{n+1}(z),
{\cal P}_j(z,w_j)=0, j=1,\cdots, n\}$ and let
$V$ be its irreducible component at $(0,0)$ which
 contains ${\Gamma_f}=\{(z,w): 
w=f(z), z (\approx 0)\in D\}$. Let $\pi$ and $\pi'$ be the natural projections
from $V$ to the first and second copies of $\bf C^{n+1}$, respectively.
Since ${\cal P}_j(0;X)=0$ has only zero solutions,
i.e, $\pi ^{-1}(0)=0$,  $\pi$  is locally proper.
Thus, we can choose a small neighborhood $U$ of $0$ such
that (a) $\pi^{-1}(U)$ is irreducible; (b)
$\pi$ is proper from $\pi^{-1}(U)$;
and (c) $\pi$ is a sheeted covering
map from $\pi^{-1}(U)\sm \pi^{-1}({\cal E})$ to $U\sm {\cal E}$
(see, for example , Theorem 1 of [Chi], pp 122).
 In the following discussion,  we make ${\cal U}_2\supset
U\supset {\cal U}_1$ and restrict $\pi$ to $\pi^{-1}(U)$.

Now, by using some standard arguments as appeared
 in the proof of Lemma 3.0 (iii)
and using the irreducibility of $\pi^{-1}(U)$,
 one can 
see that ${\cal F}(z)=\pi'\circ\pi^{-1}(z)$, for $z\in {\cal U}_2$.
(Otherwise, we can form a proper analytic subvariety of $V$ 
(of the same dimension), which is 
defined by $\{(z,w):\ w\in {\cal F}(z)\}$).
 
On the other hand, it obviously holds that ${\cal F}^{-1}(0)\cap U=
\pi(\pi'{}^{-1}(0)\cap \pi^{-1}(U))$. Thus, if $0$ is an accumulation
point of $Y={\cal F}^{-1}(0)$, then $Y\cap U$
is an analytic variety of $U$ 
 (by the elimination  theorem) and
it must have  positive dimension at $0$.
We will assume this  and seek a contradiction.

Then,   $Y\cap U$ contains
some holomorphic curve $Y^*\subset {\cal U}_1$ parametrized by
$z'=\phi(t), \zn=\psi(t)$, 
($t\in \D$, the unit disk in $\bf C^1$)
with $\phi(0)=\psi(0)=0$, $\|d\phi\|+\|d\psi\|\not =0$
for $t\not =0$.
For each $z\in Y\cap {\cal U}_1$,
 ${\cal F}(Q_z)\subset Q'_0$, by Lemma 3.1. We claim
that $\cup Q_{z}$ with $z\in Y^*$ fills in an open subset in
$\bf C^{n+1}$. This then gives us a contradiction; for we assumed that
$J_f\not \equiv 0$ and $Q_0'$ is a complex hypersurface.

To see the size of $U^*=\cup Q_{z}$ with $z\in Y^*$, 
we note that $M_1$ is also defined  by $\overline {\wt{\rho_1}(z,\overline z)}=0$.
Therefore, it can be seen that
$Q_z$  can be parametrized by
$$\zn=-\overline {\psi(t)}-\sum_{j=0}^{\infty}\overline {\wt{\phi_j}}(\overline {\phi(t)},z')
{\overline {\psi(t)}^j}
\equiv g(z',\overline {t}).$$ Here and in what follows, we use the notation $\overline {h}(z)$
for the function $\overline {h(\overline z)}$.
So, $U^*$ can be parametrized by the map 
$${\cal T}(z',\overline {t}):\ {\bf C^n}\times \D\ra {\bf C}^{n+1}$$
${\cal T}(z',\overline {t})=(z',-\overline {\psi(t)}-\sum_{j=0}^{\infty}\overline {\wt{\phi}_j}(
\overline {\phi(t)},z')\overline {\psi(t)}^j)$.
To see that ${\cal T}(z,t)$ is a biholomorphism
at certain point $(z',\overline  t)(\approx (0,0))$,
it suffices to show that ${\PP {g}\OO\PP \overline {t}}\not =0$ at 
$(z',\overline  t)$. Indeed, this can be argued as follows:
When $\psi(t)\not \equiv 0$, we can simply choose $(0',t)$ for some 
$t$ with $\psi'(t)\not =0$;
 when $\psi\equiv 0$, if ${\PP g\OO {\PP \overline  t}}\equiv 0$,
then  for any given $z'$ with $|z'|$ small,
 $\overline {\wt{\phi_0}}(\overline {\phi(t)},z')$ is independent of $t$.
 Hence it has to be  $0$; for $\phi(0)=0$ and
thus $\overline {\wt{\phi_0}}(\phi(0),z')=0$. Letting $z'=\phi(t)$, it follows that
$\overline {\wt{\phi_0}}({\overline {\phi(t)}},{\phi(t)})\equiv 0$. This contradicts the
finite type assumption on $M_1$; for it implies that $Y^*=(\phi(t),0)$
stays inside $ M_1$.
So, by making ${\cal U}_2$ small,
 it holds that ${\cal F}^{-1}(0)\cap {\cal U}_2=\{0\}$. 

Therefore, we conclude that $0$ is an isolated point of
$\pi'{}^{-1}(0)$. Hence, $\pi'$ is locally finite near $0$. Moreover,
 ${\cal F}^{-1}$, when restricted as a map from $O(0)$ to ${\cal U}_2$,
 maps a small neighborhood $0\in {\cal U}'_1$
into a small neighborhood of $0\in {\cal U}_1$.

Now,  shrinking ${\cal U}_2$ if necessary,
 we assume that $q\approx 0$, ${\cal F}^{-1}(q)\cap {\cal U}_2=
\pi(\pi'{}^{-1}(q))$ is finite. By the elimination theorem ([Chi], Theorem 1 of pp122),
 there exists an open
neighborhood $U'$ of $0$ such that $\pi'$ is proper from
$\pi'{}^{-1}(U')$ to $U'$; and away from a proper analytic set, $\pi'$
gives a sheeted covering map.
In particular, $\pi'$ is an open mapping (see [GR], Lemma 6 pp 102,
and  notice the openness follows from the onto part of that lemma).
 
In the following discussion,
we make ${\cal F}(U)\subset U'$, and we restrict $\pi'$ to $\pi'{}^{-1}(U')$.
By the above discussion,
 we also observe the fact:  $\pi^{-1}(U)\subset {\pi'} ^{-1}(U')$.

Using again the  elimination theorem
and noting that $\pi$, $\pi'$ are local analytic  covering
maps, it follows that
for any  analytic variety $ E\subset U$,
 ${\cal F}( E)\cap \pi'(\pi^{-1}(U))
= \pi'({\pi}^{-1}(E)\cap \pi'(\pi^{-1}(U)))$
is an analytic variety in  $\pi'(\pi^{-1}(U))$ with
the same dimension at the origin
(see for example, Theorem 11 E pp 68 of [Wh]; or Theorem 1, pp122, of [Chi]).
  This completes the proof of (iv).

(v)  Since $\pi$ is proper, it is also closed. Also, both are open
mappings;
 for they are local analytic covering mappings, too.
 (See [GR], Lemma 6, pp102)

For any closed subset $B$ of $U'$,
we first notice that $${\cal F}(B)\cap {U}=\pi\left ({\pi'} ^{-1}(B)\cap\pi^{-1}(U)\right).$$
Since ${\pi'}^{-1}(B)$ is closed in ${\pi'} ^{-1}(B)\cap {\pi'} ^{-1}(U')$,
it is closed in  ${\pi'} ^{-1}(B)\cap \pi ^{-1}(U)$ by the above arrangement.
By using the closeness of $\pi$, it follows that
${\cal F}^{-1}(B)$ is closed in $U$.  In particular, we see that
$\wt{M}_1\cap U={\cal F}^{-1}(M_2\cap U')$ is closed in $U$. Clearly,
${\cal F}(\wt{M_1})\subset M_2$ by part (ii) of this lemma.

Now, let $U^0=U\sm {\wt M}_1$. Then it is open. We notice
that $\PP U^0$ is contained in $\PP U \cup {\wt M}_1$. Since
 ${\cal F}=\pi'\circ\pi^{-1}$, the openness 
of $\pi'$ implies the openness of ${\cal F}$ from $U$ as
 a multiple-valued map. Therefore it sends interior points to
interior points. On the other hand,
 using the weak continuity of ${\cal F}$
(Lemma 3.0 (ii)), one sees that ${\cal F}$ maps the closure
of $U^0$ to the closure of ${\cal F}(U^0)$. Hence, a simple topological
argument shows that
 the boundary of ${\cal F}(U^0)$ is contained in
$ {\cal F}(\PP U)\cup M_2$. Since $ {\cal F}^{-1}(0)\cap {\cal U}_2
=\{0\}$, it  follows easily
that
$0$ is not contained in the closure of ${\cal F}(\PP U)$. Hence, we see that
there is a small ball, denoted by $B_0$, centered at the origin such that
$B_0\cap {\cal F}(\PP U)=\emptyset$.
Let $U^*$ be a small neighborhood of $0$ such that ${\cal F}(U^*)\subset\subset
B_0$. If $U^0\cap U^* =\emptyset$, then Part (ii) of this lemma shows that
${\cal F}(U^*)\subset {\cal F}(\wt {M_1})\subset M_2$ and thus we are done.
So, without loss of generality, we assume that $U^*\cap U^0\not =\emptyset$.
Let $D^*$ be a connected component of $ U^0$ with $U^*\cap U^0\not
 =\emptyset$. 
Note that $\PP {\cal F}(D^*)\cap B_0=\emptyset$.
Then a simple topological argument  indicates that  
 either
 ${\cal F}(D^*)\supset B_0\cap\O^+$  or
${\cal F}(D^*)\supset \O^-\cap B_0$ (see, for example, [BHR2]).
Without loss of generality, we assume the first case. Then, there exists
a point $p_0\in D^*$ and certain $\wt{f}(p_0)\in {\cal F}(p_0)$ such that
$\wt{f}(p_0)\in {\H n}^+$ and $\wt{f}(p_0)$ is close to $0$.\
 By the first part of this lemma and
Lemma 3.1, we see
that ${\cal F}(p_0)\subset \O^+$. Now, for any $z\in D^*\sm {\cal E}$,
 let $\g$ be a curve
in $D^*$ with $\g(0)=p_0$, $\g(1)=z$, and $\g((0,1))\cap {\cal E}=\emptyset$.
Then for any $\wt{f}(z)\in {\cal F}(z)$, there is a branch $\H f$
of ${\cal F}$,
 which is continuous on $\g$ and $\H f(z)$ coincides with 
$\wt{f}(z)$.  Now, since $\H{f}\circ\g$ does not meet $M_2$ by
our choice of $U^0$, and
$\H{f}\circ\g(t)\in \O^+$ for  $t$ close enough to $0$
 (this can be seen by the fact
that all limit points of $\H{f}(\g (t))$ ($t\ra 0)$
 are in ${\cal F}(p_0)$), we conclude
that $\H f(z)$ and thus $\wt{f}(z)$  have to be in $\O^+$. For
$z\in {\cal E}\cap D^*$, by passing to a limit and noting
${\cal F}(z)\cap M_2=\emptyset$, we can also see that ${\cal F}(z)\subset
\O^+$.
 Thus the proof of the last part
 of Lemma 3.3 is complete if we make ${\cal U}_1$ small. $\endpf$
\bigbreak
{\bf Remark 3.4}: The following observations will be used in the later discussions:

(a)   From the proof of Lemma 3.3 (iv), we also
notice the following property of ${\cal F}$:
There exists an open neighborhood ${\cal U}'_0$ of $0$ such that
${\cal U}_0={\cal F}^{-1}({\cal U}'_0)\cap {\cal U}_2$ is an open subset
of ${\cal U}_1$. Moreover,
 ${\cal F}|_{{\cal U}_1\ra {\cal F}({\cal U}_1)}$
 and ${\cal F}^{-1}|_{{\cal U}_0'\ra {\cal U}_0}$
 are
 multiple-valued
open mappings. 

(b)  
 Let  $(V,\pi,\pi(V))$ be the m-sheeted analytic covering space, as in
Lemma 3.3 (iv). Let $z_0\in {\cal E}\cap \pi(V)$
be close to the origin. Then, it is easy to see that in case
 $\# \pi^{-1}(z_0)<m$,
 then there are two distinct
sequences $\xi_j,
\eta_j\in \pi^{-1}(z_j)\subset V$ such that 
$z_j\ra z_0$, $\xi_j,\eta_j\ra w_0\in \pi^{-1}(z_0)$.

(c)
Still let $V$ as above.
Write $(Y,\sigma,V)$ for the standard normalization
of $V$ (see [Wh], Chapter 8). We remark that after making $V$ small,
$\sigma^{-1}(0)$ is a single point. 
Write 
 ${\cal E}_0=\pi\circ\sigma (S)$, where 
$$S=\{x\in Y: \ \h{either x is  singular or
x is smooth but} \  \h{d}_x\sigma \ \h{is singular}\}.$$ 
Write $S_0=(\pi\circ\sigma)^{-1}(\pi\circ\sigma)(S)$.
Then $\pi\circ\sigma$ is a local  biholomorphic mapping
from $Y\sm S_0$ to its image and gives a finitely sheeted covering mapping
away from the singular set.
 (see [Wh], Chapter 8
and  [GR],
pp 108). Therefore, if $V$ is singular at $0$, i.e, $f$ is not
 holomorphic at $0$, then
${\cal E}_0=\pi\circ\sigma(S)$ is an analytic variety of codimension 1
at $0$.
( This follows
from the fact that any  finitely sheeted analytic
covering space with `genuine' singular set of 
 codimension $\ge$ 2
is trivial) (see the following observation in Remark 3.4 (d)). Now,
 for each $z\in {\cal E}_0\cap O(0) (\subset {\cal E})$ and 
any small neighborhood
$U_z$ of $z$, the basic property for the normal space ([Wh], in particular
 Lemma 2F pp 254 ) indicates
that  $V$ can not be smooth at $\pi^{-1}(z)$. Therefore, $\#\pi^{-1}(z)<m
$.

(d) More generally,
let $P(z,X)=X^N+\sum_{j<N}a(z)X^j=0$ be an irreducible
polynomial equation in $X$ for $z$ near $0$.
 Let $Z^*$ be its branching locus. A point $z_0$
is called
a genuine branching point of $P=0$ if
 $P(z,X)$ can not be factorized into linear functions in $X$ in the
 Noetherian ring ${\cal O}_{z_0}$.
 Write $Z$ for the collection
of all genuine branching locus of $P$. The normalization argument as in
  $(c)$ shows that
$Z$ is an analytic variety near the origin. Moreover,
it can be seen that the following monodromy
theorem holds: Let $\g_1, \g_2$ be  loops in $O(0)\sm Z$, which are based at
$p_0\in O(0)\sm Z^*$ and are homotopic to each other  in
$O(0)\sm Z$, relative to the base point.
 Then, 
 any analytic continuation of the roots of $P$
along $\g_1$, $\g_2$, starting with the same initial value, will end
 up with the same
value when coming back to $p_0$. From this fact, if follows that in case
$Z\not=\emptyset$, $Z$ must be of codimension 1 everywhere.

\section {A local Hopf lemma and a
preservation principle of ${\cal F}$}

Before proceeding further, we need to strengthen Lemma 3.3 (iii) to
the following version: 
\bigbreak
{\bf Lemma 4.1}: $f_{n+1}(z)=z_{n+1}g(z)$ with $g(0)\not = 0$.
\bigbreak
This  sort of the Hopf lemma was established
in [BR2] in the case when the map is assumed to be smooth. 
 Since we now  have a nice control of the branches
of ${\cal F}$ and we know that
  $f_{n+1}$
is  holomorphic,  Lemma 4.1 can  be proved 
by  using the same approach and ideas as 
 in [BR2]. Because our emphasis of the present paper is the two 
dimensional case,  for completeness,
we will present a detailed proof of Lemma 4.1   when $n=1$
and leave the easy modification   in  general dimensions to  Remark 4.1$'$.

We also mention that when $f$ is a CR-homeomorphism, the proof of Lemma 4.1
readily follows from Lemma 3.3 (iii). To see this, 
write $f^{-1}=(\cdots, f^{-1}_{n+1})$. Then the same 
argument as above shows that $f^{-1}_{n+1}(w)=w^mg^*(w)$ for
some $m\ge 1$. Thus $\zn={\zn}^{km}O(1)$. So $k=m=1$.
\bigbreak
{\it Proof of Lemma 4.1}:
 We will use the approach appeared in [BR2].
We assume $n=1$.   Seeking a contradiction, we 
suppose that $f_{2}(z)=z_{2}^kg(z)$ with
$k\ge 2$.

We start with the function $G({\cal F},\l)$, which
 has the following property:
 $$\overline {f_2(\o)}=
G({\cal F}(z),\overline {\wt f_1(\o)})$$ for any $z\in Q_{\o}$ 
and $(\wt{f_1},f_2)\in {\cal F}$.
  Write the defining equation of $M_1$  in the 
following  form (see (2.1)): 
$t=\phi(z_1,\overline {z_1},s)$ with $t=2\h{Re}z_2$,
 $s=\h{Im}z_2$, and  $\phi=\rho_1-t$. 
Clearly, $\phi(0,\overline {z_1},s)=\phi(z_1,0,s)\equiv 0$. 
 Then $z\in Q_{\o}$ reads as
$${z_2+\overline {\o_2}}=\phi(z_1,\overline {\o_1},{z_2-\overline {\o_2}\over 2i}).$$
 We let $\tau=-is+{1\over 2}\phi(z_1,\xi,s)$ 
and $R(z_1,\xi,\tau)
=is+{1\over 2}\phi(z_1,\xi,s)$. As in [BR2],
using the implicit function theorem,
we can find a holomorphic function   $\mu(z_1,\xi)$ near $0$ 
 such that $R(z_1,\xi,\mu(z_1,\xi))\equiv 0$ and 
$\mu(z_1,0)=\mu(0,\xi)\equiv 0$. It is easy to verify that
$\mu(z_1,\xi)\not\equiv 0$ and thus
 ${\PP \mu\OO \PP z_1}(z_1,\xi)\not\equiv 0$; for, otherwise, it implies that
$\phi_0\equiv 0$ and thus contradicts the finite type assumption of $M_1$.
Now, one can directly verify that 
$(z_1,R(z_1,\xi,\tau))\in Q_{\overline {(\xi,\tau)}}$. So, we have
$\overline {(\xi,\tau)}\in Q_{(z_1,R(z_1,\xi,\tau))}$ and thus
$$\overline {f_2(z_1,R(z_1,\xi,\tau))}=G({\cal F}(\overline {\xi},\overline {\tau}),\overline {{\wt f_1}( z_1,R(z_1,\xi,\tau)}).$$
As in the previous section, we write $\overline {\wt f_1}(z)=\overline {{\wt f_1}(\overline {z})}$ and use the same notation for $\overline {f_2}(z)$. Then
$$\overline {f_2}(\overline {z_1},\overline {R(z_1,\xi,\tau)})=G({\cal F}
(\overline \xi,\overline \tau),\overline {\wt f_1}(\overline  {z_1},
\overline {
R(z_1,\xi,\tau)}).\leqno (4.1)$$

We still let $\cal E$
 be the branching locus of $\cal F$, i.e, ${\cal E}=\{z\in {\cal U}_3:
\ {\cal P}_1(z,X)={\PP {\cal P}_1\OO \PP X}(z,X)=0\ \h{for some}\ X\}$;
and write ${\cal F}(z)=({\cal F}_1(z), f_2(z))$.
 
 We will assume momentarily the fact 
 that $0$ is an isolated point
 of ${\cal E}\cap Q_0={\cal E}\cap \{z: \ z_2=0\}$, which will be
proved in  Lemma 5.1 (a) of $\S 5$. (In n-dimensions, we observe
the last statement of Lemma 5.1 (a) holds without Condition S).

Since $\{{\cal F}^{-1}(w)\}$ is finite for $w\approx 0$
 (Lemma 3.3 (iv)),
we see that  no branch of ${\cal F}$ is constant on an open subset
of $Q_0$. We also notice, by Lemma 3.0 (iii), that
${\cal F}_1(z_1, 0)$ consists of exactly the solutions of
 ${\cal P}_1((z_1,0);X)=0$. Using the Puiseaux expansion, we can write
${\cal F}_1(z_1,0)=g(z_1^{1/k})$, where  $g$ is  a holomorphic function
with $g(z_1)\not =0$ for each $z_1(\not =0)$ sufficiently close to $0$.

Now, for each non zero $z_1$ close to the origin,
 there is a small neighborhood $O((z_1,0))$
 near $(z_1,0)$ such that
  we can stratify $\cal F$ into several non-constant holomorphic branches. 
Next, we choose $\xi$ sufficiently small and then let $\tau$ be
 sufficiently
close to $ \mu(z_1,\xi)$.
Then $R(z_1,\xi,\tau)\approx 0$, too. After letting
$\wt{f}$ be a holomorphic branch of ${\cal F}$ over $O((z_1,0))$,
  we then can take the derivative with respect
 to $\overline  z_1$ in $(4.1)$,
to obtain
$${{\PP \overline {f_2}}\over \PP{z_1}}+{{\PP \overline {f_2}}\over 
\PP {z_2}}\overline {R'_{z_1}(z_1,\xi,\tau)}={\PP G\over \PP{\l}}({\PP \overline {
\wt {f_1}}\over \PP {z_1}}+
{{\PP \overline {\wt f_1}}\over \PP {z_2}}  \overline {R'_{z_1}(z_1,
\xi,\tau)}).\leqno (4.2)$$
Let $\tau=\mu(z_1,\xi)$ in $(4.2)$
 and notice that $f_2=z_2^{k}g$ with $k\ge 2$. 
We obtain 
$${\PP G\over \PP {\l}}({{\PP \overline {\wt f_1}}\over \PP {z_1}}(\overline {z_1},0)+
{{\PP \overline {\wt f_1}}\over \PP {z_2}}(z_1,0)  \overline {R'_{z_1}(z_1,
\xi,\mu(z_1,\xi))})=0. \leqno (4.3)$$

 By  the above discussion, 
we see that
${\PP \overline {\wt f_1}\over \PP {z_1}}(\overline {z_1},0)\not =0$ for $z_1\not = 0$.
So, for each   $z_1$ with $|z_1|$ small,
 when $\xi$ is chosen so that $|\xi|$ is sufficiently
 smaller than $|z_1|$,
 using
the fact that $R'_{z_1}(z_1,0,\tau)=0$, we see that
 $R'_{z_1}(z_1,\xi,\mu(z_1,\xi))\approx 0$. Moreover,
$${{\PP \overline {\wt f_1}}\over \PP {z_1}}(\overline {z_1},0)+
{{\PP \overline {\wt f_1}}\over \PP {z_2}}(z_1,0)
  \overline {R'_{z_1}(z_1,\xi,\mu(z_1,\xi))}\not
=0.$$
Hence it follows from $(4.3)$, that
$${\PP G\over \PP {\l}}({\cal F}(\overline {\xi},
\overline {\mu (z_1,\xi)}),\overline {\wt f_1}(\overline {z_1},0))\equiv 0.$$

On the other hand, similar to an argument in [BR2],
letting $\tau=\mu(z_1,\xi)$ in $(4.1)$, we have
$$G({\cal F}(\overline {\xi},\overline {\mu (z_1,\xi)}),\overline {\wt f_1}(\overline {z_1},0))\equiv 0.$$
Now, let $q_1=G'_{\l}(w_1,0,\l)$
and $q_0=G(w_1,0,\l)$. Notice that
 $G(w_1,w_2,\l)=-w_2+q_0-\sum_{j>0}\wt{\psi}_j(w_1,\l)w_2^j$
and $q_1=(q_0)'_{\l}$.
We see that the equations $G(w_1,w_2,\l)=0$, 
$G'_{\l}(w_1,w_2,\l)=0$ can be solved as
$w_2=q_0h^*(q_0,w_1,\l))$  and 
$q_1=q_0\psi^*(w_1,\l,q_0)$ for certain holomorphic functions $h^*$, $\psi^*$.
In the following discussions, we always let $\tau=\mu (z_1,\xi)$.

Let $w_1(z_1,\xi)=\H{f_1}(\overline {\xi},\overline {\tau})$ with
$\H{f_1}(\overline {\xi},\overline {\tau})\in {\cal F}(\overline {\xi},\overline {\tau})$,
 and let $\l(z_1)
=\overline {\wt f_1}(\overline {z_1},0)$ in the
above formulas, we obtain the following equality
$$q_1(
\H{f_1}(\overline {\xi},\overline {\tau}), \overline {\wt f_1}(\overline {z_1},0))=
q_0(\H{f_1}(\overline {\xi},\overline {\tau}),\overline {\wt f_1}(\overline {z_1},0))
\times\psi^*(\H{f_1}(\overline {\xi},\overline {\tau}),\overline {\wt f_1}(\overline {z_1},0),q_0).\leqno (4.4)$$ 

Here for clarity, we summarize 
 the situations in which $(4.4)$ makes sense and holds:

(a) $z_1$ is an arbitrarily given non zero complex number near the origin.
(b) ${\wt f}_1$ is some holomorphic branch of ${\cal F}$ over
some small neighborhood  $O((z_1,0))$ of $(z_1,0)$, whose size depends
on $z_1$.
(c) $\xi$ is taken in a small neighborhood $O(z_1,\wt{f}_1)$ of $0$, whose
size depends on the choice of $O((z_1,0))$ and the choice of $\wt{f}_1$.
(d) $\H{f}(\overline {\xi},\overline {\tau})$ is any point in ${\cal F}_1(\overline {\xi},\overline {\tau})$.

 Since $f_2\not\equiv 0$ and ${\PP \mu\OO\PP z_1}
\not
\equiv 0$,  $Z^*=\{(z_1,\xi):\ f_2(\overline \xi,\overline {\mu(z_1,\xi)})
=0\}$ is an 
analytic variety of dimension 1.
 So, we can assume that $Z^*$ has only one
component of the form $\{(a,\xi)\}$ near the origin, i.e, the one with
$a=0$.
 We therefore see that 
  $f_2(\overline \xi,\overline {\mu(z_1,\xi)})$, when regarded as a function in $\overline \xi$,
 is not  constant for each fixed $z_1\not =0$; for if it were constant,
then it would be identically $0$.
Meanwhile, noticing that 
$(\H{f_1}({\overline \xi},{\overline \tau}),{f_2}(\overline {\xi},\overline {\tau}))\in
 Q'_{(\wt{f_1}( z_1,0),0)}$,
we see that 
$$f_2(\overline \xi,\overline {\mu(z_1,\xi)})=-\wt{\psi_0}(\H{f_1}(\overline \xi,\overline {\mu(z_1,\xi)}),
\overline {\wt {f_1 }}(\overline {z_1},0)).\leqno (4.5)$$

Let $Z^{**}\subset {\bf C^2}$ be the zero set of 
$q_1(w_1,\l)-q_0(w_1,\l)\psi^*(w_1,\l,q_0(w_1,\l))$. 
If it is not of dimension 2, then it has only finitely many irreducible 
one dimensional components near the origin. 
Thus, we can assume that for any  $\l(\not =0)$ with $|\l|$ small,
$Z^{**}\cap \{(w_1,\l)\}$ 
is of dimension $0$. Returning to (4.4), one sees that
for each  $z_1\not =0$ close to the origin, $\H{f_1}(\overline \xi,
\overline {\mu(z_1,\xi)})$ takes only finitely many values
 for $\xi\in O(z_1,\wt{f})$. This obviously
 contradicts (4.5); for $f_2(\overline \xi,\overline \mu(z_1,\xi))$ can not take only
finitely many values when $\xi\in O(z_1,\wt{f}_1)$ by the above argument.

Therefore, we showed that
$(q_0)'_{\l}=q_1(w_1,\l)\equiv q_0(w_1,\l)\psi^*(w_1,\l,q_0(w_1,\l))$.
Now, applying the same argument as  in [BR2], we can see,
 by the basic theory of ODE and
by the initial condition $q_0(w_1,0)=G(w_1,0,0)=0$, that $q_0\equiv 0$.
This contradicts the finite type assumption of $M_1$. $\endpf$

\bigbreak
{\bf Remark 4.1$'$}:  The
 above proof can be adapted to the general dimensions as in [BR2],
 almost without any
change. However, we would like to
include here some details just for 
 convenience of the reader:

We first replace $z_1 $ by $z'\in {\bf C}^n$ and define, in the same manner, 
the functions $\tau$, $\mu$, $R$. We then have, by Lemma 3.1, the following

 $$\overline {f_{n+1}(z',R(z',\xi,\tau))}=
G({\cal F}(\overline {\xi},\overline {\tau}),\overline {\wt {f^*}(z',R(z',\xi,\tau)});$$
for $\overline {(\xi,\tau)}\in Q_{(z',R(z',\xi,\tau))}$ with $\xi\in {\bf C^n}$.

Let ${\cal E}'={\cal E}_0\cap Q_0$. By Lemma 5.1 (a), it is a proper
analytic variety of $Q_0$. In the same manner,
we let $z'\in Q_0\sm {\cal E}'$ be close to $0$ and choose
a small neighborhood $O((z',0))$ of $(z',0)$ such that
${\cal F}$ can be stratified into several (non-constant)
holomorphic branches over $O((z',0))$. Arguing in the same way,
we have for each $l\le n$,
$${{\PP \overline {f_{n+1}}}\over \PP{z_l}}+{{\PP \overline {f_{n+1}}}\over 
\PP {z_{n+1}}}\overline {R'_{z_l}(z',\xi,\tau)}
=\sum_{j=1,\cdots,n}{\PP G\over \PP{\l_j}}({\PP \overline {
\wt {f_j}}\over \PP {z_l}}+
{{\PP \overline {\wt f_j}}\over \PP {z_{n+1}}}  \overline {R'_{z_l}(z',
\xi,\tau)}).$$

Let $\tau=\mu(z',\xi)$. Assume that $f_{n+1}=z_{n+1}^kg$ with $k>1$.
One then has, for $\xi$ with $\|\xi\|$ sufficiently small, that
$$\sum_{j=1,\cdots,n}{\PP G\over \PP{\l_j}}({\PP \overline {
\wt {f_j}}\over \PP {z_l}}(\overline {z'},0)+
{{\PP \overline {\wt f_j}}\over \PP {z_{n+1}}}(z',0)  \overline {R'_{z_l}(z',
\xi,\tau)})=0.$$

Since ${\cal F}^{-1}$ is a finite-to-finite map,
$f^*(z',0)$ is a finite-to-one map from $O(z',0)\cap Q_0$ too.
 Moving to a nearby point if necessary,
we can assume that $$\h{det}(
{\PP \overline {\wt f_j}\over \PP {z_l}}(\overline {z'},0))_{1\le j,l\le n}\not =0.$$
Reasoning  as in the proof of the two dimensional case and making
$|\xi|$ sufficiently small,
we have for each $l<n+1$
$${\PP G\over \PP {\l_l}}({\cal F}(\overline {\xi},
\overline {\mu (z',\xi)}),\overline {\wt f^*}(\overline {z'},0))\equiv 0,$$
and
$$G({\cal F}(\overline {\xi},\overline {\mu (z',\xi)}),\overline {\wt f^*}(\overline {z'},0))\equiv 0.$$
Here, $\|\xi\|$ is sufficiently small.
Let $q_l=G'_{\l_l}(w',0,\l)$ and
$q_0=G(w',0,\l)$.  Similarly, we   have
$$q_l(w',\l)=q_0(w',\l)h^*_l(w',\l,q_0(w',\l)),$$
where $w'=\H{f}^*(\overline {\xi},\overline {\tau})$,
$\tau=\mu(z',\xi)$,
 $\l(z')=\overline {\wt{f}^*}(\overline {z'},0)$, and $h_l^*$ is holomorphic in
its variables.
(See the notation  introduced in the proof of Lemma 4.1).

Now, a similar argument as in the two dimensional case (involving the use of 
the n-dimensional version of (4.5)) indicates
that
$$q_l(w',\l)\equiv q_0(w',\l)h_l^*(w',\l,q_0)\leqno(4.6)$$ for
$(w',\l)\approx 0$ and each $l<n+1$.
 Indeed, this can also be seen as follows:
Let $z'_0\in {\bf C^n}$ be chosen as above
such that $\l(z')=\overline {\wt{f}^*}(\overline {z'},0)$ is an open mapping
from $O((z'_0,0))\cap Q_0$ to $\bf C^n$. Fix each $z'\in O((z'_0,0))\cap Q_0$
and consider the map $H$ which sends each $\xi(\approx 0)\in {\bf C}^n$ to 
$(\H{f}^*(\overline {\xi},\overline {\mu(z',\xi)}),f_{n+1}(\overline {\xi},\overline {\mu(z',\xi)}))$.
 As argued before, we see that
$H(\xi)\in Q'_{({\wt{f}^*}({z'},0),0)}$ and is a finite to one map, too.
Now, let $\pi^*$ be the natural projection from $\bf C^{n+1}$ to its first
n-copies of $\bf C$. Then, it follows easily 
that $\pi^*\circ H$ is an open mapping.
Hence, from the uniqueness theorem of holomorphic functions, we conclude that
$(4.6)$ holds identically. Finally, the same argument as in [BR2]
gives a contradiction.

\bigbreak
{\bf Lemma 4.2}: After shrinking ${\cal U}_1$, it holds that
${\cal F}({\cal U}_1\cap D)\subset \O$ and 
${\cal F}({\cal U}_1\cap D^c)\subset \O^c$. Here, $\O$ and $\O^c$ stay in
different sides  of $M_2$.
\bigbreak
{\it Proof of Lemma 4.2}: 
We will still use the $\l -$function $G({\cal F},\l)$ introduced in Lemma 2.1.
Let 
$$\rho^*(z,\overline w)=1/N\h{Re}\left(\sum_{\wt{f(w)}\in {\cal F}(z)}\left(
\overline {f_{n+1}(w)}-
G({\cal F}(z),\overline {\wt{f(w)}})\right)\right),$$
where $N$ is the generic counting number of ${\cal F}(z)$.
Then by Lemma 2.1 and the Riemann extension theorem, we see
 that $\rho^*(z,\overline {w})$ is holomorphic in $(z,\overline {w})$. In particular,
$\rho^*(z,\overline  z)$ is real analytic near $0$.
We claim that $\rho^*(z,\overline z)$ is $<0$, when ${\cal F}(z)\subset \O^+$;
$>0$ when ${\cal F}(z)\subset \O ^-$; and $=0$ when
 ${\cal F}(z)\subset M_2$ (see Lemma 3.3 (a) for related notations).
In fact, for each fixed $\wt{f}(z)\in {\cal F}(z)$, we have
$$\h{Re}\left (\overline {f_{n+1}(z)}-G({\wt f},\overline {\wt f(z)})\right)
=\h{Re}\wt{\rho_2}\left({\wt f}, \overline {\wt f(z)}\right)$$
$$=\rho_2\left({\wt f}, \overline {\wt f(z)}\right)
\h{Re}h_2({\wt f(z)}, \overline {\wt f(z)}))=
\rho_2({\wt f}, \overline {\wt f(z)}))(1+o(\|\wt{f(z)}\|).$$
Thus, the claim follows from Lemma 3.3 (v).

Next, applying  Lemma 4.1, one sees that $$\rho^*=2\R (f_{n+1}(z))+o(|z|)=
\R (g(0)z_{n+1})+o(|z|)$$ with $g(0)\not = 0$.
  Thus,  $\h{d}_0\rho^*\not =0$. Therefore,
$\rho^*$  serves as a real analytic
 defining function of $M_1$ near $0$.
 Without loss of generality,
let us assume that $\rho^*(z,\overline z)<0$ for $z\in D$. Then $\rho(z,\overline z)>0$ for
$z\in D^c$. Now, by making ${\cal U}_1$ small, the above argument
shows that ${\cal F}(D\cap {\cal U}_1)\subset \O^+$ 
and ${\cal F}(D^c\cap {\cal U}_1)
\subset \O^-$.
So, we can let the $\O$ in Lemma 4.2 be $\O^+$.
 Of course, $\O^c$ is then $\O^{-}$.   $\endpf$
\bigbreak
{\bf Remark 4.3} (a). In what  follows, we will assume,
without loss of generality, 
 that $\O=\O^+$. Then
we  observe that the above results give the fact $g(0)>0$. To see this,
we  write
$f_{n+1}(z)=\a z_{n+1}+o(\|z\|)$
with
$\a=\a_1+\sqrt{-1}\a_2\not =0$.
 Letting  $t=x+\sqrt{-1}kx$,
then $\rho^*(t)=(\a_1-\a_2k)x+o(|x|)$. Notice that
$\rho^*(x)<0$ for any given real number $k$ and $x(\approx 0)<0$. We conclude that
$(\a_1-\a_2k)>0$. Thus, we see that $\a_2=0$ and $\a_1>0$.

(b). Another immediate application of Lemma 4.2 is that for certain small
 neighborhoods $U^*$ and ${U^*}'$ of $0$, $f$ is proper from $U^*\cap D$
to ${U^*}'\cap \O$. In fact, let $B$ be a sufficiently small ball centered at $0$.
By certain well-known construction (see [BC], for
example) and Lemma 3.3 (iv), 
we can simply take $U^*$ to be the connected component of the set
$\{z\in B\cap D:\ 
f(z)\not\in \PP {B\cap D}\}$ such that $f(U^*)$ contains $\O\cap O(0)$
with $O(0)$ sufficiently small.

\section{Completion of the proof of Theorem 1.1}

We now proceed to the proof of our main technical theorem. For simplicity,
we retain all the notation which we have set up in the previous sections.

 As before, denote by ${\cal A}$ 
the  map, which sends:
$w\in {\cal U}_1$ to the finite set
 $A_{w}$. $w_0$ is called a {\it separable} point of ${\cal A}$ if
$A_{w_0}=\{w_j\}_{j=1}^{N^*}$ satisfies the following property:
There exist open neighborhoods ${ O}({w_j})$ of $w_j$ ($j=1,\cdots, N^*$)
such that ${\cal A}(w)\cap {O}_{w_j}=\{w\}$ for
any $w\in { O}_{w_j}$.
 We write ${\cal W}=\{
w\in {\cal U}_1:\ w \ \h{is not a separable point of}\ {\cal A}\}$.
We say that $M_1$ satisfies {\it Condition S} at $0$, if a holomorphic
change of coordinates can be chosen so that $M_1$ is still defined by
an equation of the form as in (2.1) and (2.1)$'$; moreover, the following holds:

\noindent (a) $\cal W$ is contained in an analytic variety $E$ spread over
$(z_2,\cdots,z_{n+1})$-space, i.e, $E$ is defined by an equation of the
form: $z_1^{N^*}+\sum_{j=0}^{N^*-1}a_j(z_2,\cdots,z_{n+1})z_1^{j}=0$ with
$a_j(z_2,\cdots,z_{n+1})$ holomorphic near the origin.

\noindent (b) For each small $t>0$, there is a sufficiently small positive number $\epsilon(t)$
such that for each $b^*=(b_2,\cdots,b_n)$ with $|b^*|<\epsilon(t)$,
the set ${E}\cap \{(z_1,b^*,-t): z_1\in {\bf C}^1\}$ is contained in some connected component of
  $D\cap \{(z_1,b^*,-t): z_1\in {\bf C}^1\}$.

Similarly, we can define ${\cal A}'$ and ${\cal W}'$ and speak of {\it Condition S} for $M_2$. In this case, we use $\O$ to replace the role of $D$
(see Remark 4.3(a)).

We will see in the following lemma that ${\cal W}$ and ${\cal W}'$
 can  be used to control
the branching locus of ${\cal F}$.

\bigbreak

{\bf Lemma 5.1}: 
\noindent (a) Suppose that $M_1$ and $M_2$ satisfy Condition S at the origin. Then  
${\cal F}({\cal E}_0\cap O(0))\subset  {\cal W}'$,
${\cal E}_0\cap O(0)\subset {E}$, and
  ${\cal E}_0\cap O(0)\cap \{(z_1,\cdots,z_n, 0)\}$ is a proper variety near
the origin of the
 $\{(z_1,\cdots,z_n,0)\}$-space. Here $E$ is chosen as in the definition
 of Condition S, and ${\cal E}_0$ is as defined in Remark 3.4 (c).
(When $n=1$, we always have ${\cal F}({\cal E})\subset {\cal W}'$ and
thus $0$ is an isolated point of ${\cal E}\cap Q_0$).
\noindent (b) $M_1$ and $M_2$ satisfy Condition S,
 when $n=1$ or when $M_1$ and $M_2$ are
rigid at the origin.

\bigbreak
{\it Proof of Lemma 5.1}:  
We first prove Part (a) of the lemma:

Let $z\in {\cal E}_0$ be sufficiently close to $0$. Then, by the
definition of ${\cal E}_0$ and Remarks 3.4 (b) (c), it follows that there is
a sequence $z_j\ra z$ such that
we can find two
sequences $\{\eta_j\}$ and $\{\xi_j\}$ with $\xi_j,\eta_j\in {\cal F}(z_j)$
$\eta_j\not=\xi_j$, but $\xi_j,\eta_j\ra w_0$ for some $w_0$.
By Lemma 3.1, it then follows that $\xi_j\in A'_{\eta_j}$. Hence,
$w_0$ is not a separable point of ${\cal A}'$. Again by Lemma 3.1, 
we conclude that ${\cal F}(z)\subset {\cal W}'$. {\it Notice that for
this part, no Condition S is necessary}.  

We next prove that ${\cal E}_0\cap O(0)\subset {E}$. 
To this aim, we assume, without loss of generality,
that ${\cal E}_0\not=\emptyset$.
  Let
$\cup_j{\cal E}_j$ be the  union of the irreducible components of
${\cal E}_0$ near $0$. We recall by Remark 3.4 (d), that 
each ${\cal E}_j$ then must be of codimension 1 at $0$.
 Now, since ${\cal F}({\cal E}_j)
\subset {\cal W}'$ is also an analytic variety of codimension 1 at $0$
(Lemma 3.3 (iv)),
the assumption of Condition S for $M_2$ 
 indicates that ${\cal F}({\cal E}_j)\cap \O\not 
=\emptyset$ near $0$. Therefore, it follows from Lemma 4.2, that
${\cal E}_j\cap D$ is not empty in any neighborhood of $0$. Now, pick
any point $z\in {\cal E}_j\cap D$, which is sufficiently close to $0$.
 Then, as above, we can find a sequence
 $z_j\ra z$ and
 $\{\eta_j\}$,  $\{\xi_j\}$ with $\xi_j,\eta_j\in {\cal F}(z_j)$
$\eta_j\not=\xi_j$, but $\xi_j,\eta_j\ra w_0\in \O$. By making
$z$ sufficiently close to $0$, we can assume by Remark (4.3)(b), that
 $f(z_0)=w_0$ for some $z_0$ close to $0$. Moreover, $f$ is proper
from a small neighborhood $U$ of $z_0$ and $\overline {U}\cap f^{-1}(w_0)=\{z_0\}$.
 Since $f$ is single valued,
 there are two distinct points $a_j$, $b_j\in U$ such that
$f(a_j)=\xi_j$ and $f(b_j)=\eta_j$. Now, the following
 Claim 5.2 indicates that
$a_j\in A_{b_j}$. Since $a_j,b_j$ clearly converge to $z_0$
as $j\ra \infty$, we conclude that
$z_0\in {\cal W}$. Again by Lemma 3.1 and Claim 5.2, we see that
$z_0,z\in {\cal W}$. By the arbitrariness of the choice of $z$ and
the assumption of Condition S,
 this fact implies that ${\cal E}_j$ has an open subset contained
in ${\cal W}$ and thus in $E$.
 From the uniqueness of analytic varieties, it follows that
${\cal E}_j\subset { E}$ for any $j$.
\bigbreak
{\bf Claim 5.2}: Let $M_1$, $M_2$, and $f$ be as before.
We have that
$b\in A_a$ if $Q'_{{\cal F}(a)}=Q'_{{\cal F}(b)}$, where
 $a,b\in {\bf C}^{n+1}$
 are
sufficiently close to $0$. (For this claim, no Condition $S$ is required). 
\bigbreak
{\it Proof of Claim 5.2}:
 Let $a$ and $b$ be as in the claim.
We then need 
to show that $Q_a\cap {\cal U}_2=Q_b\cap {\cal U}_2$.
 To this aim, we consider the function $\wt{\rho_2}({\cal F}(z),
\overline {\wt{f}(\o)})$, for any given $\wt{f}(\o)\in {\cal F}(\o)$.  By 
Lemma 2.1, $\wt{\rho_2}({\cal F}(z),\overline {\wt{f}(\o)})$ is well defined and
 holomorphic
in $(z,\l)$ with $\l=\overline {\wt{f}(\o)}$ for $(z,\l)\in O_{z}(0)\times O_{\l}(0)$. 
 Meanwhile $\wt{\rho_2}({\cal F}(z),\overline {\wt f(\o)})
=cz_{n+1}+c\overline {\o}_{n+1}+O(|z||\l|+o(\o))$ for some $c>0$, by Remark 4.3 (a).
 Thus,
 for each  $\o$ sufficiently close to $0$,
 since $|\wt{f}(\o)|\approx 0$,  the implicit function theorem indicates
 that
 $\wt{\rho_2}({\cal F}(z),\overline {\wt f(\o)})=0$ defines a connected
 complex hypersurface 
 $W_{\o}$
  in some small  neighborhood $U^{\#}$ of $0$, whose size
is independent of the  parameter $\o$ with $|\o|<<1$.
 Notice, by Lemma 3.1, that
 ${\cal F}(Q_{\o})\subset Q'_{\wt{f}(\o)}$ near $0$.
We see that 
 $Q_{\o}\cap U^{\#}\subset W_{\o}$. 
Thus it follows that $W_{\o}=Q_{\o}$ in 
 $ U^{\#}$ when $|\o|<<1$; for
$Q_{\o}$ is also a connected complex submanifold of codimension 1 in
${\cal U}_2$.
 
Now, if $Q'_{{\cal F}(b)}= Q'_{{\cal F}(a)}$,
 we then
 conclude that
  $$\wt{\rho_2}({\cal F}(z),\overline {{\cal F}(b)})=0$$ and ${\wt \rho_2}({\cal F}(z),
\overline {{\cal F}(a)})=0$ define the same variety. 
 Hence, when $a,b\approx 0$, we  conclude, by Lemma 3.0 (iv),
 that $Q_{a}=Q_{b}$; for  both of them have a small open subset
of $W_{a}\cap U^{\#}=W_{b}\cap U^{\#}$ in common.
 This gives
 the proof of Claim 5.2. (We mention that for the proof Claim 5.2,
no assumption of Condition S is required).
\bigbreak

The last statement in (a) follows clearly
  from the above  facts   and
the assumption of Condition $S$ for the hypersurfaces. More precisely,
we notice that ${\cal F}({\cal E}_0\cap O(0))\subset {\cal W}'$ and
$ {\cal F}(Q_0)\subset Q'_0$.  From Lemma 3.3 (iv)
and the fact that ${\cal W}'\cap \{(z_1,\cdots,z_n,0)\}$ is contained
in a proper subvariety
(of codimension at least 1) in $Q'_0$, it follows that ${\cal E}_0\cap
 Q_0$ is also of codimension at least $1$ in $Q_0$.
{\it Actually, using the same proof as in Part (b) for the rigid case, and
using a simple projection lemma, we notice that this part does not require
Condition S neither}.  Also, we notice that the argument for this part does not
use
  Claim 5.2, i.e, the use of Lemma 4.1, neither.

\bigbreak
Now, we proceed to the proof of (b):

We first  consider the case of $n=1$.
As before, we let $M_2$ be defined by 
$$z_{2}+\overline {z_{2}}+\sum_{j=0}^{\infty}\wt{\psi_j}(z_1,\overline {z_1})z_{2}^j=0;
\h{or}\ 
z_{2}+\overline {z_{2}}+\sum_{j=0}^{\infty}\overline {\wt{\psi_j}}(\overline {z_1},{z_1})
\overline {z_{2}^j}=0.$$
For each $b=(b_1,b_{2})\approx 0$, $A'_b=\{(w_1,w_{2}):\ Q'_{w}\cap
 {\cal U}'_2=Q'_{b}\cap {\cal U}'_2\}$.
As did in Lemma 3.3, by letting $z_1=0$ in the defining equations of
$Q'_w$ and $Q'_b$, one sees that $w_{2}=b_{2}$. Now,
$Q_{w}$ can also be defined by
$$z_2=-\overline { b_2}-\sum_{j=0}^{\infty}\overline {{\wt \psi_j}}(\overline  w_1,z_1)\overline {b_2}^j,$$
where the notation $\overline {\wt \psi}$ is the same as explained before.
Write $$\overline {\wt \psi_j}(\overline {w_1},z_1)=
\sum_{\a>0}\overline {\Xi_{j\a}}(\overline {w_1})z_1^{\a}.$$
Then $z_2=-\overline {b_2}-\sum_{\a}\sum_{j}\overline {\Xi_{j\a}}(\overline  w_1)\overline {b_2}^j z_1^{\a}$.
So, the equation: $Q_{w}'=Q'_{b}$ can be written as
$$\sum_{j}\overline {\Xi_{j\a}}(\overline  w_1)\overline {b_2}^j=\sum_{j}\overline {\Xi_{j\a}}(\overline  b_1)
\overline {b_2}^j$$ or
$$\sum_{j}{\Xi_{j\a}}( w_1){b_2}^j=\sum_{j}{\Xi_{j\a}}( b_1){b_2}^j,$$
for any $\a$.

Choose $k_0$ to be the smallest integer so that for certain $\a_0$, 
 $\Xi_{0\a_0}(w_1)=a^*w_1^{k_0}+o(|w_1|^{k_0})$ with $a^*\not =0$.

 Write $P(w_1,b_2)=\Xi_{0,\a_0}(w_1)+\sum_{j>0}\Xi_{j,\a_0}(w_1)b_2^j$.
For $(b_1,b_2)$ close to the origin,
we define 
  $$E=\{(b_1,b_2):\ \h {for some}\ w_1^0,\ \ 
P(w_1^0,b_2)=P(b_1,b_2),\ {\PP P\over \PP w_1}(w_1^0,b_2)=0\}.$$
Then we first claim that
near the origin,
${\cal W}\subset E$.
For this, let $w=(\xi,b_2)\not \in E$, and let
$C(w)=\{(z_1,b_2):\ P(z_1,b_2)=P(\xi,b_2)\}$.
 Then by the definition of $P$,
one has $A'_w\subset C(w)$. Write $A'_w=\{(\tau_j,b_2)\}_j$ and suppose that
$w\in W_0$. Then for some $j^*$, by the definition of $\cal W$,
 there are two sequences $\{(\eta_j, c_j)\}$
and $\{(\xi_j,c_j)\}$ such that (a) $\eta_j\not =\xi_j$; (b) $\eta_j, \xi_j
\ra \tau_{j^*}$ for some index $j^*$, $c_j\ra b_2$;
 and (c) $A'_{(\eta_j, c_j)}=A'_{(\xi_j,c_j)}$. 

From (c) and the definition of $P$, it follows that
$P(\eta_j,c_j)=P(\xi_j,c_j)$. Moreover,  by the Taylor expansion:
$P(\eta_j,c_j)=P(\xi_j,c_j)+{\PP P\OO \PP w_1}(\xi_j,c_j)(\eta_j -\xi_j)
+o(|\eta_j-\xi_j|)$; we obtain 
${\PP P\OO \PP w_1}(\xi_j,c_j)=o(1)$. Letting $j\ra 0$, we  therefore
 conclude that
${\PP P\OO \PP w_1}(\tau_{j^*},b_2)=0$, which together with the
 fact $P(\tau_{j^*},b_2)=P(\xi,b_2)$ implies that  $w\in E$.
 This is a contradiction
and thus proves our claim.

 We next claim that $E$ is an analytic variety
spread over $b_2$-axis. To see that, we let 
$$E^{*}=\{(b_1,b_2,w_1)\in O(0):
P(w_1,b_2)=P(b_1,b_2),\  {\PP P\over \PP w_1}(w_1,b_2)=0\}.
$$
 Then $E$ is the projection of $E^{*}$ to the first two copies of $\bf
C^1$. This projection map is obviously finite to one near the origin and thus
is locally proper. So, $E$ is an analytic variety  near 0.
Meanwhile, for each $b_2$, from
${\PP P(w_1,b_2)\over \PP w_1}=0$, we  can solve out only
 finitely many $w_1$'s with $|w_1|$ small.
Moreover, for any solution $w_1$,
 we can easily see that $|w_1|\ale |b_2|^{{1\over k_0-1}}$.
Thus it follows  that there are also only finitely many $b_1$'s so that
$P(w_1,b_2)=P(b_1,b_2)$ with
${\PP P(w_1,b_2)\over \PP w_1}=0$. Furthermore, 
one can easily see that
   $|b_1|\ale |b_2|^{1/k_0}$.
 Therefore, we see that $E$ is
spread over $b_2$-axis.
 From the Weiestrass preparation
theorem, it follows that $E$
 can be defined near $0$ by an equation with the form  as given in
Part (a) of Condition S.
 Meanwhile, it is also easy to see that
$E\cap ({\bf C^1}\times \{b_2\})\subset \D_{|b_2|^{\a^*}}\times
\{b_2\}$ for $|b_2|$ sufficiently small,
 where $\D_{r}$ is used to stand for the disk in $\bf C^1$ with
center at the origin and with radius $r$,  $\a^*$ is a constant
between $1\over k_{0}+1$ and $1\over k_0$. By our choice of $k_0$,
it follows easily for sufficiently small $|b_2|$, that  
$ \D_{|b_2|^{\a^*}}\times
\{b_2\}\subset D$ when $b_2<0$ , and
 $\D_{|b_2|^{\a^*}}\times
\{b_2\}\subset D^c$ when $b_2>0$.
This completes the proof of Part (a) in case $n=1$.

%
Finally, we show that $M_1$ satisfies Condition S when $M_1$ is rigid.
For this purpose, we assume that
$M_1$ is defined by an equation of the form:
$z_{n+1}+\overline {z_{n+1}}+\sum_{|\a|,|\b|>0}a_{\a\b}z'{}^{\a}\overline {z'{}^{\b}}=0$.
Let $P_{\a}(z')=\sum_{\b}\overline {a_{\a\b}}z'{}^{\b}$. Then
for $a=(a',a_{n+1})$ $b=(b',b_{n+1})$ close to the origin,
$Q_{a}=Q_{b}$ if and only if $a_{n+1}=b_{n+1}$ and
$P_{\a}(a')=P_{\a}(b')$. By the finite type assumption, it follows that
the common zero of $P_{\a}'$s is $0$ near the origin. Using the
Nullstellantz theorem, we see that for some finitely many indices
$\{\a_j\}_{j=1}^{m}$, the locus of $\{P_{\a_j}\}_{j=1}^{m}$ is
 also zero. Now, for a small neighborhood $U$ of 0, let 
$U'=U\cap \{(z_1,\cdots,z_n,0)\}$. Define
$\Lambda$ from a small neighborhood $U'$ of $0'$ to $\bf C^{m}$ by
$\Lambda(z')=(P_{\a_j}(z'))_{j=1}^{m}$. 
 Then $\Lambda$ is finite to $0$ and proper from $U'$,
after suitably shrinking of $U$. So, by the Remmert theorem, we conclude
that the set $E'=\Lambda^{-1}\left(\L(\{z\in U: \ \h{d} {\Lambda}\
 \h{does not have maximal rank at} \
z\})\right)$ is a proper variety of $U'$ (see [DF1] for
a similar argument).
 It is easy to see that ${\cal W}\subset
E$, with $E=(E'\times {\bf C^1})\cap {\cal U}_1$.
Now, after a linear change of coordinates in the 
$(z_1,\cdots,z_{n})$-space, we see that $E$ is contained in an
analytic variety defined by an equation of the form:
$z_1^{N^*}+\sum_{j<N^*}c_j(z_2,\cdots,z_n)z^j_1=0$ with
$c_j(0)=0$. Since this equation 
does not contain $z_{n+1}$-variable, 
it follows easily that $M_1$ satisfies Condition $S$ at the origin.

\bigbreak

We now give the following  lemma, which  indicates how
the information on the branching locus can be used to obtain
the holomorphic extendibility. We  use $\D_{a}$ to denote
the disc in $\bf C^1$ centered at $0$  with radius $a$.
\bigbreak
{\bf Lemma 5.3}: Let $M$ be a hypersurface
 near the origin of $\bf C^{n+1}$ defined by an equation
of the form $\rho(z,\overline  z)
=z_{n+1}+\overline {z_{n+1}}+\phi(z,\overline {z})=0$ with $\phi(z,z)=o(|z|)$.
 Let $D$ be the  side defined by $\rho(z,\overline {z})<0$.
 Assume  that $g$ is a holomorphic
function over $D$ that admits a
holomorphic correspondence extension to $\D_{\d_0}\times\D^{n}_{\d}$
for
some small $\d_0, \d>0$,
i.e,
there exist holomorphic functions $a_j\in \h{Hol}(\D_{\d_0}\times
\D_{\d}^{n})$
such that $P(z,g(z))=g^N(z)+\sum_{j=0}^{N-1}a_jg^j(z)=0$
for $z\in D$.
 Here we assume that $P(z,X)$ is irreducible in $X$,
 and $\D_{\d_0}\times\D_{\d}^{n}\cap D\subset\subset D$.
 Denote by $Z^*$   the branching locus
of $g$ and write $Z$ for the genuine branching locus of
$g$ as defined in Remark 3.4 (d). Assume that $Z$ is defined by
 $\wt{P}(z_2,\cdots,z_{n+1};z_1)=z_1^{N_1}+\sum_{j=0}^{N_1-1} d_j(z_2,
\cdots,z_{n+1})z_1^j=0$
for certain $d_j(z_2,\cdots,z_{n+1})\in \h{Hol}(\D ^{n}_{\d})$
with $d_j(0)=0$.
 Also, suppose that  for each small positive number $t$, there is
a positive number $\epsilon(t)$ such that when $|b_j|<\epsilon(t)$
with $j=2,\cdots, n$,
 $Z\cap D^*(b^*,t)$ is contained in some
  connected open subset of $D\cap D^*(b^*,t)$,
 where $D^*(b^*,t)=\{z:\ z=(z_1,b_2,\cdots,b_n,-t): z_1\in \bf C^1\}$.
 Then $g$ admits a holomorphic
 extension near $0$.

\bigbreak
{\it Proof of Lemma 5.3}:
Without loss of generality, we assume that $$\wt{P}(z_2,\cdots,z_{n+1};X)=0$$
 has no double roots,
i.e, $\wt{P}(z_2,\cdots,z_{n+1};X)$
 and ${\PP \wt{P}\OO \PP X}$ are relatively prime.  Denote
 by $\wt{Z}$ the branching locus of $\wt{P}$, which  is a subvariety of
 $ \D_{\d}^{n}$. Now, we choose a small $t>0$ and
$b^*=(b_2,\cdots,b_{n})$ with $|b_j|<\epsilon(t)$ such that
$(b^*,-t)\not\in \wt{Z}$
and $Z^*\cap {\bf C^1}\times \{(b^*,-t)\}$ is discrete.
 Write $B_{\eta}(b^*,-t)$ for the {\it closed} ball in $\bf C^{n}$
 centered at
$(b^*,-t)$ with radius $\eta$. Then, when $\eta$ is sufficiently small,
$B_{\eta}(b^*,-t)\cap \wt{Z}=\emptyset$ and $0\not\in B_{\eta}(b^*,-t)$.
 Construct a  pseudoconvex
domain $\wt{B}(t)\subset\subset
\D_{\d}^{n}$ such that the following holds:
 (i) $0\in \wt{B}(t)$, $\wt{B}(t)\supset B_{\eta}(b^*,-t)$; and
(ii) there is a contraction map $T(\cdot,\tau): \wt{B}(t)\ra \wt{B}(t)$
($\tau\in [0,1]$) such that
 $T(\cdot,0)=\h{id}$, $T(\wt{B}(t),1)=
 B_{\eta}(b^*,-t)$,
and for any $\tau$,
 $T(\cdot,\tau)$ is  the identical map when restricted to 
 $B_{\eta}(b^*,-t)$.

Now, using the hypothesis on $Z$ and shrinking $\d$ if
necessary, we can find $ R>0$ with $R<\d _0$ such that
 $ (\D_{R}\sm \D_{{1\over 2}R})\times \D_{\delta}^n\cap Z=\emptyset$.
Write $\O_0=\left(\D_{R}\times\wt{B}(t)\right)\sm \overline {
\left(\D_{{1\over 2}R}\times (\wt{B}(t)\sm B_{\eta}(b^*,-t))\right)}$, which can
also be written as $$\left(\D_{R}\times \h{int}(B_{\eta}(b^*,-t))\right)
\cup \left((\D_{R}\sm \overline {\D_{{1\over 2}R}})\times
 \wt{B}(t)\right).$$ Here, for a subset $K$, we write $\h{int}(K)$ for
the collection of the interior points of $K$. 
Then, it is easy to see that any (closed) loop
in $\O_0\sm Z$  based  at certain $z_0\in
\left (\D_{R}\times \h{int}(B_{\eta}(b^*,-t))\right )\sm Z^*$
  can be deformed, relative to the base point and without cutting $Z$,
 to a loop in
$\left(\D_{R}\times
 B_{\eta}(b^*,-t)\right)\sm Z$. 
Indeed, let $\g=(\g_1,\g^*)$ be  a loop in $\O_0\sm Z$ based at $z_0$.
 Then
the above described deformation can be performed by
$\left(\g_1(t),T(\g^*(t),\tau)\right)$.

 We next set up the following notation:
 For a given connected subset $D'$ (near the origin) with $D'\sm Z$ 
pathwise-connected, a point $z_0\in D'\sm Z^*$,
 and a complex number $w_0$ with $P(z_0,w_0)=0$, we define
$\h{Val}(z_0,w_0;D')$ to be the collection of $w'$s, which
is obtained as follows:
There are a loop $\g\subset D'\sm Z,$ based at $z_0$, and
a continuous section $\cal G$ from $[0,1]$
to $({\cal O},\pi,{\bf C})$ with $P(z;{\cal G})\equiv 0$
such that each representation of
 ${\cal G}$(0) takes value $w_0$ and each
 representation
of ${\cal G}(1)$ takes value $w$ at $z_0$.
We  write ${\cal P}(z)=\{w:\ P(z,w)=0\}$.
When $\V(z_0,w_0;D')$ is independent of the choice of $w_0$, we simply use the
notation $\V(z_0;D')$. In particular, this is the
 case when $\#\V(z_0,w_0;D')=N$.
\bigbreak
{\bf Claim 5.4}:
  $\#\V(z_0; \D_{R}\times B_{\eta}(b^*,-t))=N$, i.e,
 $$\V(z_0; \D_{R}\times(B_{\eta}(b^*,-t))=\{\h{all
 solutions of}\ P(z_0,X)=0\}.$$
 \bigbreak
{\it Proof of Claim 5.4}: Let $\O_0$ be as defined above.
 We then notice that the holomorphic hull of $\O_0$ is $\D_{R}
\times \wt{B(t)}$. In fact, for any $\phi\in\h{Hol}(\O_0)$, 
 the following Cauchy integral gives the holomorphic extension of $\phi$
to $\D_{R}\times \wt{B}(t)$:
$$
{1\over 2\pi i}\int_{|\xi|=\wt R}{\phi(\xi,z_2,\cdots,z_{n+1})
\over \xi-z_1}\h{d}\xi,\ 
\h{with}\ \wt R\ra R.$$
Now, we fix a point $\ z_0\in\O_0\sm Z^*$
 and a number $w_0$ with $P(z_0,w_0)=0$.
For $z\in \O_0\sm Z^*$,  let $\g\subset \O_0\sm Z$ be
 a curve connecting $z_0$ to $z$.
Denote by $w$ be the value of the branch of ${\cal P}$ at $z$
which is continuous along $\g$ and takes
the initial value $w_0$.
   Write
 $\V(z,w;\O_0)=\{w_1,\cdots,w_{k(z,w)}\}$.
 Then, by the monodromy theorem (see Remark 3.4 (d)),
one sees that $k(z,w)$ is  constant, say $k'$, for all $z, w$. Also,
$\V(z,w;\O_0)$ is independent of the choice of $\gamma$, though $w$ does.
 Now,
$P_0(z,X)=\prod_{j=1}^{k'}(X-w_j(z))$ will give a well-defined polynomial
 in $X$ with coefficients holomorphic and bounded in $\O_0\sm Z^*$.
 Applying the Riemann extension theorem,
 we conclude that $P_0(z,X)$ is a polynomial in $X$ with coefficients
 holomorphic in $\O_0$. The above observation then indicates
that the coefficients of $P_0$ are holomorphic in $\D_{R}\times
 \wt{B}(t)$, too. 
So, by the irreducibility of $P$, we conclude that $k'=N$.
 
Now, from the argument preceding this claim and the monodromy theorem, 
 the proof of Claim 5.4 follows. $\endpf$
\bigbreak
{\bf Claim 5.5}:
Let $b^*$ and $t$ be as chosen above and write $p_0=(b^*,-t)$. Let
$\epsilon(t)$ be sufficiently small. Let
$z_0=(z_1,p_0)\in \D_{R}\times \{p_0\}\sm Z^*$ for a certain $z_1$.
Then 
$\#\V(z_0;\D_{R}\times\{p_0\})=N$.
\bigbreak
We first mention that the proof of Lemma 5.3 follows easily from
Claim 5.5. Indeed,
by the hypothesis, the finite set
$Z\cap \{\D_R\times \{p_0\}\}$ stays
 in some {\it connected} open subset of
$D\cap \{{\bf C^1}\times \{p_0\}\}$, say $D^*$.
Since $\# (Z\cap \{{\bf C^1} \times \{p_0\}\})$ is finite,
 we can choose a simply connected smooth domain $D^{**}\subset\subset D^*$
with $Z\cap ({\bf C^1}\times \{p_0\})\subset D^{**}\cap ({\bf C^1}\times
\{p_0\})$
(see Remark 5.6 (a) below).
  Therefore, by the  Schoenflies theorem,
there is a retraction from  $\D_R\times \{p_0\}$  to $\overline {D^{**}}$
that deforms $\D_{R}\times \{p_0\}\sm D^{**}$ to the boundary of $D^{**}$. Thus
any loop in $\D_{R}$ that is based at $z^*\in D^{**}\sm Z^*$ can be
 deformed, relative to the base point, to $\overline {D^{**}}$
 without cutting $Z\cap {\D_R}\times \{p_0\}$.
 Hence, by the monodromy theorem, it follows that $\#\V({z^*};D^*)=N$
for any given point $ z^*\in D^{**}$. But this is impossible; for $f$ is
 holomorphic in a small neighborhood of $D^*$. 
So, to conclude the proof of Lemma 5.3,
 it suffices for us to prove Claim 5.5.

We now turn to the proof of Claim 5.5.
By our choice of $B_{\eta}(b^*,-t)$, we know  that
$$\wt{P}(z_2,\cdots,z_{n+1};z_1)=z_1^{N_1}+\sum _{j<N_1}d_j(z_2,\cdots,z_{n+1})
z_1^j=0$$ has roots 
$$\{\phi_1(z_2,\cdots,z_{n+1}),\cdots,\phi_m(z_2,\cdots,z_{n+1})\},$$
 which can be arranged such that
   $\phi_j(z^*)\in \h{Hol}(B_{\eta}(b^*,-t))$
 and $\phi_j(z^*)\not =\phi_l(z^*)$ for $j\not =l$.
 So, for a given $z^*\in B_{\eta}(b^*,-t)$,
$(\D_{R}\times \{z^*\})\cap Z=\{(\phi_j(z^*),z^*)\}_{j=1}^{m}$.
 Let $p_j^0=\phi_j(p_0)$. We will show that there exists a 
homeomorphism  $\pi_{z^*}$:
$\D_{R}\ra\D_{R}$,
depending continuously on $z^*\in B_{\eta}(b^*,-t)$,
 so that $\pi_{z^*}(\phi_j(z^*))=p_j^0$. If this can be done, then any loop 
$\g(t)=(\g_1(t),\g^*(t))$, which is based at $z_0$ and stays
 in $(\D_{R}\times  B_{\eta}(b^*,-t))\sm Z$
 can be 
 deformed to a closed curve  in $\D_{R}\times \{p_0\}$
  by the following map:
$$\b(t,\tau)=\left(\pi^{-1}_{\tau p_0+
(1-\tau)\g^*(t)}\circ\pi_{\g^*(t)}(\g_1(t)),
 \tau p_0+(1-\tau) \g^*(t)\right)$$
with $\b(t,0)=\g(t)$, $\b(t,1)\in \D_{R}\times \{p_0\}$, and 
$\b(0,\tau)=\b(1,\tau)=z_0$. Meanwhile, from the construction,
 one can verify that $\{\b(t,\tau)\}\cap Z=\emptyset$. 
 Thus, from the monodromy theorem stated in Remark 3.4 (d),
it follows that $\#\V (z_0;\D_R\times \{p_0\})=N$.

To  construct the above mentioned $\pi_{z^*}:\D_{R}\ra\D_R$.
We notice that $$\min_{j\not = l,z^*\in B_{\eta}(b^*,-t)}
\{|p_j(z^*)-p_l(z^*)|\}>0.$$ By using a M\"obius transform,
 we  can assume that $p_1(p_0)=p_1^0=0$. 
Suppose that for an index $j<m$, we can find a homeomorphic self map
$\sigma_j(\cdot, z^*)$ of $\D_{R}$,
which depends continuously on $z^*$ and 
maps $p_l(z^*)$ to $p_l^0$ for $l\le j$. 
For the brevity of the notations, let us still write
$p_j(z^*)$ for the transformed points: $\sigma_j(p_j(z^*),z^*)$. 
 We then wish to find
a  homeomorphic self map  $\wt{\sigma_{j+1}}(\cdot,z^*)$ of $\D_R$,
 depending continuously
on $z^*\in B_{\eta}(b^*,-t)$,
 such that $\wt{\sigma_{j+1}}(p_l^0,z^*)=p_l^0$ for $l\le j$
and 
$\wt{\sigma_{j+1}}(p_{j+1}(z^*),z^*)=p^0_{j+1}.$
 If this can be done, letting $
{\sigma_{j+1}}(\cdot,z^*)=\wt{\sigma_{j+1}}(\cdot,z^*)\circ 
\sigma_j(\cdot,z^*)$,
and  applying the induction argument,
 we see the existence of the aforementioned
map $\pi_{z^*}$. 
Indeed,
 the existence of $\wt{\sigma_{j+1}}(\cdot,z^*)$
clearly follows from the following Remark 5.6 (b) and by making
$\epsilon(t)$ sufficiently small:
\bigbreak

{\bf Remark 5.6}: In this remark, we would like to
 say a few words about two  simple facts
in the set topology that were used in the proof of Lemma 5.2:

(a) Let $D$ be a connected domain in $\bf C^1$.
Then, for any given finitely many points $\{p_j\}_{j=1}^{k}\subset D$,
 there is a smoothly
bounded
simply connected domain $D^*\subset\subset D$ such that $\{p_j\}_{j=1}^{k}
\subset D^*$.

Indeed, we can find a curve $\g$ in $D$ connecting all of the $p_j$'s,
which,
 by the piecewise linear approximation, can be assumed to be piecewise linear.
Making  further division to the image of $\g$ if necessary,
 we can assume that $C=\g ([0,1])$ is a    linear graph (or, a 
finite CW-complex of dimension 1),
which contains $\{p_j\}$ in the set of its vertices.
Now, a well-known result in Topology indicates that we can find a
 maximal
tree $\wt C$ inside $C$, which, by definition, is simply connected
and contains all the vertices
 of $C$.
 Furthermore, the tubular neighborhood basis
theorem shows the existence of a small neighborhood $\wt D(\subset  D^*$)
 of $\wt C$, that can be retracted to  ${\wt C}$
 and thus is simply connected, too.
(See for example,
 J. F. Hudson: {\it Piecewise Linear Topology}, W. A. Benjamin, Inc, 1969;
in particular, Theorem 2.11, pp57 ).  Now, we have a 
Riemann mapping $\sigma$ from $\wt D$ to $\D$.
 Then obviously, $\sigma^{-1}(\D_{r})$
 does our job when $r<1$ is sufficiently close to $1$. $\endpf$.
\bigbreak 
 (b)
 Let $\{a,a_1,\cdots, a_k\}$ be $k+1$ 
distinct points in $\D$. Suppose that $p(z)$ is a continuous map 
from the ball the closed unit ball $ B_n\subset {\bf C^n}$ to $\D$ 
such that $\epsilon=\min_{j=1,\cdots,k;z\in B_n}|p(z)-a_j|$ is a sufficiently
small positive number,
and $|p(z)-p(z_0)|<{1\over 4}\epsilon$. 
Then there exists a homeomorphic self map 
$\sigma(\cdot,z)$ of $\D_{R}$, 
which depends continuously on $z$, such that $\sigma(p(z),z)=a$ and
$\sigma(a_j,z)=a_j$ for $j=1,\cdots, k$.
\bigbreak
The proof of this fact can be seen as follows:
Pick a continuous retraction $\pi$ from
$B_n$ to $z_0$.
 We will then construct a homeomorphic self map $\sigma(\cdot,z)$ of $\D$,
that depends continuously on $z\in B_n$, such that the following holds:
(a) 
$\sigma(a_j,z)=a_j$ $(j=1,\cdots, k)$; (b)   
$\sigma(\cdot,z_0)=\h{id}$; (c) 
$\sigma(p(z),z)=p(z_0)$ for $z\in B_n$.
To this aim, 
we assume that
 $p(z_{0})=0$ to simplify the notation. 
Notice that
 $|p(z)-p({z_0})|<\epsilon/4$ and thus
$p(z)$ is contained in the disk $\D_{\epsilon/4}$. We also observe, by the
hypothesis, that $a_j\not\in \D_{\epsilon}$. Now, we let $\tau(\cdot,z)$
 be the
M\"obius transform from $\D_{\epsilon/4}$ to itself that maps $p(z)$ for
$z$ to $p(z_0)$. Meanwhile, we impose the condition that
$\tau'(p(z),z)>0$ to make it unique. From the explicit formula of $\tau$,
it is clear that $\tau(\cdot,z)$ depends continuously on $z\in B_n$.

For $w=\epsilon/4 e^{\sqrt{-1}\theta}$ with $0\le \theta<2\pi$, write
$\tau(w,z)={\epsilon\OO 4}e^{\sqrt{-1}r(\theta,z)}$,
 where $r(0,z)\in [0,2\pi)$, 
$r(2\pi,z)-r(0,z)=2\pi$, and $r(\theta,z)$ is a strictly increasing function 
of $\theta$ for fixed $z$. Moreover, 
$r(\theta,z)$ depend continuously on $(\theta,z)$. 
The homeomorphism $\sigma(\cdot,z)$ from $\D$ to itself
can be defined as follows:

(i)  $\sigma(\cdot,z_0)=\h{id}$. 
For $z\in B_n$,
we let $\sigma(w,z)=w$, for $w\not\in \D_{\epsilon}$; 
$\sigma(w,z)=\tau(w,z)$, for $w\in \D_{\epsilon/4}$;
and otherwise
$$\sigma(w,z)= |w|\h{exp}\left({\sqrt{-1}{4\OO 3\epsilon}\left(
(|w|-\epsilon/4)\theta+(\epsilon-|w|)h(\theta,z)\right)}\right),$$
where $w=|w|\h{e}^{\sqrt{-1}\theta}$ with $|w|\in [\epsilon/4,\epsilon]$.

It can be easily 
verified that $\sigma$ possesses all the properties we imposed. 
\bigbreak

{\it Completion of the proof of Theorem 1.1}.
We now are ready to apply the previously
 established results to present a proof of Theorem 1.1.
 Indeed Theorem 1.1 follows from the
following Theorem 1.1$'$ and Lemma 5.1 (b).
\bigbreak
{\bf Theorem 1.1$'$}: Let $M_1\subset \bf C^{n+1}$ and $M_2\subset \bf C^{n+1}$ be real analytic
hypersurfaces of finite D-type. Suppose that $f$ is a CR mapping
from $M_1$ to $M_2$. Let $p\in M_2$.
 Assume that $f$ extends as a holomorphic correspondence near
$p$ and assume furthermore that $M_1$, $M_2$ satisfy Condition S at $p$, and $f(p)$, with respect to $D$ and $\O$ as introduced before, respectively.
Then $f$ admits a holomorphic extension near $p$.
\bigbreak
{\it Proof of Theorem 1.1}$'$:
For each component $f_j$ of $f$, we see that the genuine
branching locus of the
polynomial equation defining $f_j$, denoted by  ${\cal E}_0^{(j)}$,
 is contained in
$
{\cal E}_0$. 
By Lemma 5.1, it then follows that for sufficiently small $t>0$,
one can always find $\epsilon(t)>0$ such that when
 $|b_j|<\epsilon(t)$ ($j=2,\cdots, n)$,
${\cal E}^{(j)}_0\cap \{(z_1,b_2,\cdots, b_n,-t)\}$ is contained in some
connected component of $D\cap  \{(z_1,b_2,\cdots, b_n,-t)\}$. Now, applying
Lemma 5.3, we conclude that $f_j$ extends holomorphically
 across $0$ for each $j$.
$\endpf$
\bigbreak
{\bf Remark 5.7}: By examing the proof,
 it can be seen that Theorem $1.1'$ 
also holds when $M_1$ and $M_2$ are merely assumed to  be  essentially 
finite
at the origin in the sense that $A_0=A'_{0}=\{0\}$ (see [BJT] or [DW]).

\section {Proofs of Theorem 1.2  and Theorem 1.4}

Now, we pass to the proofs of  Theorem 1.2 and Theorem 1.4. 
We first let $M_1$, $M_2$, and $f$ be as in Theorem 1.2. 
We also assume that $f$ is not constant.
Then we claim that $f$
must be an algebraic map. To see this, we first stratify $M_2$
 into the disjoint union of the sets
 $M_2^{\pm}$, $S_2$, $S_1$, and $S_0$. Here $M^{\pm}_2$
 is the set of points in $M_2$ where the Levi form
 is non-zero;
$S_2$ is a locally finite union of
  $2-$dimensional real analytic totally real
submanifolds 
in  $M_2\sm M_2^{\pm}$; 
$S_1$ is a locally finite union of  1-dimensional 
real analytic submanifolds in $M_2\sm (M^{\pm}\cap S_2)$;
 and $S_0$ is a locally finite subset. For the existence of such
a semi-analytic stratification, we refer the reader to a related
explanation appeared in [DFY]. 
Assume that $0\in M_1$ and $f(0)=0$. 
Since the algebraicity is a global property, we need only
to show the map is algebraic at some small piece of $M_1$.
Using the finite type assumption of $M_1$, 
we can assume the existence of some non empty open piece $U_0$
in $M_1^{\pm}$, where $M_1^{\pm}$ denotes the set of points in $M_1$
where the Levi form do not vanish. 
When $f(U_0)\cap M_2^{\pm}\not =\emptyset$, by  
results of  Pinchuk-Tsyganov [PS],  Lwey [Le], and Pinchuk [Pi], 
we see that $f$ is actually real analytic at some point $p$ in $U_0$. 
Thus applying a result of Baouendi-Rothschild [BR3],
 we see $J_f\not\equiv 0$ near $p$.
Now, the algebraicity follows from a theorem of Webster [We].

When $f(U_0)\subset S_2$, by shrinking $U_0$ if necessary, we assume that
$f$ extends to a certain side, say $D$, of $M_1$.
As in the proof of Lemma 2.1, 
we can find a totally real submanifold
$E\subset U_0$ and a wedge $W^+\subset D $ with edge $E$. Moreover,
we can assume that $f(E)$ is contained in a connected piece of $S_2$.
Now, by the reflection principle,
we conclude that $f$ admits a holomorphic extension
near $E$. Since $f$ maps a neighborhood of $E$ in $U_0$ to a two dimensional
manifold, we easily see that $J_f\equiv 0$ near $E$.
 This contradicts the  non constant assumption for $f$, by [BR3]. 
Similarly, we can exclude the case: $f(U_0)
\subset S_1\cup S_0$. Summarizing the above arguments, we can conclude
 the algebraicity
of $f$.

We  also recall a result in [DF2], which states that any proper holomorphic map
between two bounded algebraic domains in $\bf C^n$ ($n>1)$ is continuous up
to the boundary. Meanwhile, it is also algebraic by [We] (see also [DF2]).
Thus, the proofs of Theorem 1.2, Corollary 1.2$'$, and
 Theorem 1.4
now clearly follow from Theorem 1.1 and the following  general result: 
\bigbreak
{\bf Lemma 6.1}: Let $M_1$ and $M_2$ be two connected
algebraic hypersurfaces of
finite D-type in $\bf C^n$ ($n>1$) and
assume that $D$ is bounded domain with $M_1$ as part of its boundary.
 Suppose that $f$ is an algebraic holomorphic map
 from $D$ to $\bf C^n$, that is continuous up to $D\cup M_1$, and maps
$M_1$ into $M_2$.
 Then for each point $p\in M_1$,
$f$ extends as a holomorphic
 correspondence to a neighborhood of $p$.
\bigbreak
{\it Proof of Lemma 6.1}: 
 Let $M_1$ and $M_2$ be defined by two real polynomials $\rho_1(z,\overline z)$ and 
$\rho_2(z,\overline z)$, respectively.  Without loss of generality, we 
assume  that $p=0$ and $f(0)=0$.
By the hypothesis,
 we suppose that the $j^{\rm th}$ component $f_j$ of $f$ 
satisfies the irreducible polynomial equation 
$P_j(z,X)=\sum_{l=0}^{N_j}a_{jl}
(z)X^l=0$. 
As mentioned before, by a result of Baouendi-Treves [BT] (see also [Tr]), 
we can assume, without loss of generality, that any CR function defined
over 
  $M_1$ can be extended  to $D$. 
 
 Write $V$ for the irreducible 
variety in
$\CC^n\times\CC ^n$, defined by $P_j(z,w)=0$ for
$j=1,\cdots,n$, which  contains the graph $\Gamma _f$ of $f$ over $D$.
 Write
$\pi$, $\pi'$ for the nature projections from $V$ to the first and the second copies
of $\bf C^n$, respectively. 
Define E to be the collection of points
which are either the zeros of  $ a_{jN_j}(z)$ for some $j$,
or are the branching points of of $P_j(z,X)$ for certain
$j$.
Notice that $E$ is an analytic variety of codimension $\ge 1$.
\bigbreak
{\bf Claim 6.2}: Fix a point  $z_0\in (M_1\sm E)\cap {\cal U}_1$
 and choose a sufficiently
small ball $B$ centered at $z_0$. Consider the multiple-valued extension
${\cal F}_j$ of $f_j$ from $B\cap D$ to $B\cup (D^c\cap {\cal U}_1)\sm E$
with ${\cal U}_1$ as chosen before (see $\S 3$ for a precise definition
of ${\cal F}_j$) . If there is a constant $C$ such that
$|\wt{f}(z)|\le C$ for any $z\in
B\cup (D^c\cap {\cal U}_1)$ and $\wt{f}(z)\in {\cal F}_j(z)$,
 then
   $f_j$, near $0$, satisfies 
a polynomial equation with leading coefficient 1.
\bigbreak
{\it Proof of Claim 6.2}:
For each $z\in
B\cup (D^c\cap {\cal U}_1)\sm E$,
  let $${\cal F}_j(z)=\{f_j^{(1)}(z),\cdots,
f_j^{(k(z))}(z)\}$$.
First, by the monodromy theorem,
 $k(z)$ is constant.
Now, suppose that all elements in ${\cal F}_j(z)$ is bounded by a constant
independent of $z$. As we did before,  using
the Riemann extension theorem, we conclude that
 $\wt{P_j}(z,X)=\prod_{l}(X-f^{(l)}_j(z))$ is a polynomial 
in $X$ with coefficients holomorphic and bounded in $D^c\cap {\cal U}_1$.
 Now,  by the aforementioned result of Baouendi-Treves and Trepreau, 
one sees that the coefficients of $\wt{P_j}$ can be 
extended holomorphically to a neighborhood of $0$.
$\wt{P_j}$ has leading coefficient 1 and annihilates $f_j(z)$
when $z\in {\cal U}_1$ . $\endpf$
\bigbreak
Now, seeking a contradiction, suppose  that Lemma 6.1 is false.

Then Claim 6.2 indicates that
 no matter how we shrink ${\cal U}_1$ and make
$z_0$ close to $0$,
 there  exists a fixed index $j_0$ such that one can always find
 some point $z^1\in
B\cup (D^c\cap {\cal U}_1)\sm E$ with $|\H{f_{j_0}}(z^1)|>1$ for certain
${\H f}_{j_0}(z)\in {\cal F}_{j_0}(z)$. Observe that
for a given curve $\g\subset B\cup (D^c\cap {\cal U}_1)\sm E$
with $\g(0)=z_0$ and  $\g(1)=z^1$, there is a  certain
holomorphic branch $\wt{f_{j_0}}$ of ${\cal F}_{j_0}$ along $\g$, which
 takes
value $f_{j_0}(z_0)$ at $z_0$ and takes value $\H{f_{j_0}}(z^1)$
at $z^1$. 

On the other hand, by the Whitney approximation theorem
(see J. Milnor: {\it Differential Topology}, Lecture Notes by
Munkress; in particular, Theorem 1.28 and Lemma 1.29)
 and by noting
that $B\cup (D^c\cap {\cal U}_1)\sm E$ can be retracted
to $D^c\cap {\cal U}_1\sm E$, we notice
 the existence of
 a smooth simple  curve $\g^*$
  connecting
$z_0$ to $z$ with $\g((0,1])\subset 
 (D^c\cap {\cal U}_1)\sm E$ and $\g$  transversal to $M_1$ at $z_0$.
Moreover
$\g^*$ is homorpotic to $\g$ (relative to the base point)
in $B\cup (D^c\cap {\cal U}_1)\sm E$.

So, using the Monodromy theorem 
we can simply assume  that
 $\g=\g^*$. 
As what we did in Lemma 3.1,
we can thicken $\wt{g}=\g\cup {\cal R}(\g)$ and then apply
the uniqueness theorem of holomorphic functions and the invariant property of
Segre varieties
 to conclude   
that $f(Q_z\cap O(z^*))\subset Q'_{\wt{f}(z)}$, 
for any $z\in\g$, where
the branch $\wt{f}$ of $\cal F$ is determined by the initial value condition
${\wt f}(z_0)=f(z_0)$.

Since $f(0)=0$, 
for any small   $\epsilon>0$,  
by  changing the size of ${\cal U}_1$,
we see  the existence of a sequence $\wt{p}_j\ra 0$
 such that $|\wt{f}(\wt p_j)|=\epsilon$. 
In particular, we see that the cluster set  $Cl_{\wt f}(0)$ of $\wt{f}$ 
at $0$ contains points with  norm $\epsilon$.
 Here, the branch
${\wt f}$ is as explained
above.
We now are ready to use the following argument
  to show that this is impossible,
which was first used in [DF2] for a different purpose.

Let $\H{V}$ be the algebraic closure of $V$ in
 $\bf P^{n}\times \bf P^{n}$ (see [Wh]). 
 Without loss of generality, we can further assume 
 that $\H{V}$ is irreducible and contains $\Gamma_f$. 
 Let $\H{\pi}$ and $\H{\pi'}$ be the nature projection from
  $\H{V}$
  to the first and the second copies of $\bf P^n$, respectively.
Since $\H{V}$ is of dimension $n$, 
   away from a proper subvariety, we see that 
  $\H\pi$ and $\H\pi'$ are local biholomorphisms. 
  For each $z$, let 
  $Q_z
  =\{w\in {\bf C^n}\cap {\cal U}_1:\ \rho_1(w,\overline z)=0\}$ as before, 
  and let $\H{Q}_z$ be
the compactification of $Q_z$ in $\bf P^n$. 
Now, for each $w\in Cl_{ \wt{f}}(0)$, 
the cluster set of ${ {\wt f}}$ at $0$,
let $z_j\ra 0$ with $\wt{f}(z_j)= w_j\ra w_0$.  
Then $z_j\in D^c$ and
$f(Q_{z_j}\cap O(z_j^*))\subset Q'_{w_j}$. 
This gives
that $\H{\pi'}\circ\H{\pi}^{-1}(\H{Q_{z_j}})\supset Q'_{w_j}\cap {\cal U}'_2$,
where ${\cal U}'_2$ is as defined in $\S 2$. 
Now, letting $j\ra \infty$, we see that $\H{\pi'}\circ\H{\pi}^{-1}(\H{Q_{0}})
\supset Q_{w'}\cap {\cal U}'_1$. Notice that
${\pi'}\circ\H{\pi}^{-1}(\H{Q_0})$ is an analytic variety of ${\bf P}^n$ of
codimension $1$, by the Remmert mapping theorem.
 Hence, it follows that
 $\H{Q}_{w'}\cap {\cal U}'_2$ has to be contained in
 one of 
the finite irreducible components of 
 $\pi'\circ\pi^{-1}(\H{Q_0})$. Since $\#A_{w}'<\infty$ ($w\in {\cal U}'_2$), 
 for each connected
 component
of   $\pi'\circ\pi^{-1}(\H{Q_0})$, there are only finitely many $w$'s
 with $|w|$ small such that the above property holds.
 This implies that $\#( Cl_{\wt f}(0)\cap O(0))<\infty$ 
 and contradicts our assumption.
  The proof of Lemma 6.2 is complete.
 $\endpf$

\section {Extending  proper holomorphic mappings in $\bf C^2$- Proof of Theorem 1.3}

With Theorem 1.1 at our disposal,
we now  use the ideas appeared in [DFY] and, in particular, [DP] to 
 present a proof of the following Theorem 1.3 $'$, 
 which together with Theorem 1.1 gives Theorem 1.3. We would like
to mention again the whole strategy is similar to that in [DP] (see also
closely related work in [DFY]) for the study of biholomorphic maps.
\bigbreak
{\bf Theorem $1.3'$}: Let $D$ and $D'$ be two smoothly 
bounded  domains with real analytic boundaries in $\bf C^2$.
 Let $f$ be a proper holomorphic map from $D$ to $D'$.
 Suppose that $f$ extends continuously up to
 $\PP D$. Then $f$ extends as  a holomorphic correspondence across $\PP D$.
 By Theorem 1.1, it thus follows that $f$ admits
 a holomorphic extension to $\overline {D}$.
\bigbreak

To start with, we let $M$ be a real analytic hypersurface of finite type   
in $\bf C^2$,
and let $\O$ be a domain with $M$ as part of its smooth boundary. 
Then one has the following  semi-analytic stratification as introduced
in  [DFY]:
$$M=M^+\cup M^-\cup T_2^+\cup T_2^-\cup T_2^{\pm}\cup T_1\cup T_0.$$
 Here 

\noindent (a) $M^+$ ($M^-$) is the set of strongly pseudoconvex (strongly 
pseudoconcave, respectively) boundary points of $\O\cup M$;

\noindent (b) $T_2^+$ ($T^-_2$) is a locally
finite union of 2-dimensional totally real analytic submanifolds, where
the vanishing order of the Levi form of $M$,
as the boundary of $\O\cup M$,
 is even and positive (negative, respectively). 
Therefore, $\O$ is pseudoconvex near $T^+_2$ and pseudoconcave near $T_2^-$;

\noindent (c) $T_2^{\pm}$ is a locally finite union of 
 2-dimensional totally real submanifolds, where the vanishing order of the
 Levi form
 is  odd. 
 It is known (see [DF3] or [BCT]) 
 that any CR functions defined near $T^{\pm}_2$ can be holomorphically extended to both sides of $M$
 near $T_2^{\pm}$;

\noindent (d) $T_1$ is a
locally finite union of one dimensional real analytic curves; and
$T_0$ is a locally finite subset of $M$.

Now, let $D$ and $D'$ be as in Theorem 1.3$'$. After making 
the above type of the stratifications, we have the following disjoint unions:
$$D=\PP D^+\cup\PP D^-\cup T_2^+\cup T_2^-\cup T_2^{\pm}\cup T_1\cup T_0,$$
$$D'=\PP {D'}^+\cup\PP {D'}^-\cup {T'}_2^+\cup {T'}_2^-\cup {T'}_2^{\pm}\cup T'_1\cup T'_0.$$
In what follows, we write $\Sigma'=
\PP {D'}^+\cup\PP {D'}^-\cup {T'}_2^+\cup {T'}_2^-\cup {T'}_2^{\pm}$, and write
$\Sigma$ for the boundary points in $\PP D$ where $f$ extends holomorphically.
\bigbreak
{\bf Lemma 7.1}: (a) For any  boundary point $p\in \PP D$, 
if $f(p)\in \Sigma'$,
then $f$ extends holomorphically across $p$. (b) 
$f$ extends almost everywhere in $\PP D$  and
$\Sigma\cap T_2^+$ is dense in $T_2^+$. 
\bigbreak
{\it Proof of Lemma 7.1}:   We first note
that in case $q=f(p)\in \PP D'{}^{-}\cup {T'}_2^{\pm}\cup {T'}_2^-$,
then $p$ stays inside the holomorphic hull of $D$ (see [BHR2]).
Thus $f$ extends automatically at $p$.
We now assume 
that $q=f(p)\in \PP {D'}^+\cup {T'}_2^+$. 
Pick a small pseudoconvex piece $M'\subset \PP D$ of $q$, and  
let $M$ be the connected component of $f^{-1}(M')$ which contains
$p$. Then since $f(\PP D)\subset \PP D'$, it is clear that
$f^{-1}(q)\subset\subset M$.
 We notice that $M$ has to be pseudoconvex, by
the properness and boundary continuity of $f$.
 Hence, using a result of Bell-Catlin [BC], it follows
that $f$ is smooth at $p$.
By [BBR], 
we therefore conclude that $f$ extends holomorphically across $p$.

Next, if  for some point $p$ and a small
open subset $U_p\subset\PP D$ near $p$,
$f(T_2^+\cap U_p)\subset T_1'\cup T'_0$,
 then it is easy to see that
there is a wedge $W^+$ with edge
$U_p\cap T^+_2$ such that $W^+\subset D$.
Indeed, after a local change of coordinates, one can assume that
$p=0$, $T^+_2=i{\bf R}^2$, and $T_0^{(1,0)}$ is given by $z_1$-axis.
Then, one can simply take $W^+=\{(z_1=x_1+iy_1,z_2=x_2+iy_2)\approx (0,0):\
x_1<0, |x_j|>|y_j|\}$.

Notice that $T'_1\cup 
T'_0$ is contained in the union of finitely many analytic arcs and finitely many points. Using the continuity of $f$ and
 applying the reflection principle,
 one sees that $f$ extends holomorphically in a small neighborhood of
 some point
 $z\in U_p\cap T_2^+$. 
Now, notice that
$J_f\not\equiv 0$ for $w\in U_p\cap T_2^+$ ($\approx z$).
 This gives a contradiction; for
$f$ maps the two dimensional manifold $U_p\cap T^+_2$ into
a one dimensional semi-analytic subset.
Now, from the above argument, it follows that
$\Sigma\cap T^+_2$ is dense in $T^+_2$.

 Similarly,
one can see that $f$ extends holomorphically
 across almost every point
 in $\PP D$, as did in the beginning of $\S 6$.
$\endpf$
\bigbreak
Now, we fix $p_0\in\PP D$ and $q_0=f(p_0)$. We will show that
$f$ extends as a holomorphic
correspondence near $p_0$. We can also assume, without
loss of generality [BT] [Tr], that every CR function defined on $\PP D$
near $p_0$ can be extended holomorphically to the side in $D$.

For each $p\approx p_0$ and $q\approx q_0$
with $q_0=f(p_0)$,
we introduce Segre
neighborhood systems near $p$ and $q$, as did
in $\S 2$: ${\cal U}_3(p)\supset\supset {\cal U}_2(p)\supset\supset {\cal U}_1(p)
\ni p$;
${\cal U}'_3(q)\supset\supset {\cal U}'_2(q)\supset\supset
 {\cal U}'_1(q)\ni q$.
As in [DFY] and [DP] (in particular, see [DP]),
 we define the ``potential"  correspondence of $f$ 
near $p_0$ 
as follows:
$${\cal V}(p_0,q_0)=\{(w,w')\in {\cal U}_1(p_0)\times {\cal U}_1'(q_0):\ 
f(Q_w(p_0)\cap {\cal U}_2(p_0)\cap D)\supset\  _{w^{'*}}Q'_{w'}(q_0)\}.$$
Here, as in $\S 2$, we write $w'{}^*$ for the reflection point ${\cal R}'(w')$ of
$w'$; and we use 
$ _{{w'*}}Q'_{w'}$ to denote the germ of the variety $Q'_{w'}$ at $w^{'*}$.
Also, the  Segre varieties are chosen in the indicated Segre systems.

In what follows, we will drop the reference points for the Segre systems and
Segre varieties to simplify the notation,
when there is no confusion arising.

By the invariant property of Segre surfaces, 
it can be seen that if $f$ extends
holomorphically in a small  neighborhood $G$ of some point 
$p\in {\cal U}_1(p_0)\cap \PP D$, then
${\cal V}(p_0,q_0)\supset \{(w, w')\approx (p,f(p)):\ w'=f(w),\ w\in G\sm{\overline  D}\}$. 

The key point in proving Theorem 1.3$'$
 is to show that ${\cal V}(p_0,q_0)$ extends as an 
analytic variety in ${\cal U}_1(p_0)\times{\cal U}_1'(q_0)$
 with the first projection map $\pi: 
 {\cal V}(p_0,q_0)\rightarrow {\cal U}_1(p_0)$
 locally proper. For this purpose, we will use the strategy of 
Diederich-Pinchuk [DP]
 of  applying the extension lemma of 
 Bishop ( [Chi]).

Before arguing further, we need recall some useful facts:

\noindent (7.a)
 $$\PP D\sm({\cal U}_2(p_0)\cap \PP D\cap Q_T)$$ is dense in $\PP D\cap
 {\cal U}_1(p_0)$, 
 where
$Q_T=\cup_{z\in (T_1\cup T_0)\cap {\cal U}_1(p_0)}Q_z(p_0).$

\noindent (7.b) ([DP]): Making ${\cal U}_1'$ small,
then $(T'_1\cup T'_0)\cap {\cal U}'_1(q_0)$ is contained in 
 a pluripolar set $\wt{T}$ in ${\cal U}'_1$. That is, there is
a plurisubharmonic function $h\ (\not\equiv -\infty)$ such that
$T'_0\cup T'_0\subset {\wt T}=\{z\in {\cal U}'_1:\ h(z)=-\infty\}$. 
(This  follows from the fact that the semi-analytic subset 
$(T'_1\cup T'_0)\cap {\cal U}'_1$ 
is contained in a one dimensional real analytic subset [BM], which can be 
complexified.)

\noindent (7.c) Fix ${\cal U}_1'$ . Using the properness and the boundary
continuity
of $f$, one can always choose ${\cal U}_2$ small enough so that $\cal V$
has no limit points in $({\cal U}_1\sm   {\overline D})\times \PP {\cal U}_1'$.

\noindent (7.d) ([DP]) For each $p\in T_2^+$,
 after making ${\cal U}_2(p)$ small,
it holds that $Q_p\cap {\cal U}_2(p)\cap\PP D\subset \PP D^+\cup \{p\}$.

\noindent (7.e) For each $p\in \PP D^+$, after making ${\cal U}_2(p)$ 
small and then making ${\cal U}_1(p)$ sufficiently small,
it then holds that
 $Q_{\wt p}\cap {\cal U}_2(p)\cap\PP D=\{\wt p\}$ for any
 ${\wt p}\in {\cal U}_1(p)$. 

\noindent (7.f) ([BR2]) Let $M_1$ and $M_2$ be two connected oriented
real
analytic hypersurfaces in $\bf C^n$ of finite D-type.
Let $f$ be a non constant
 holomorphic map from $M_1$ into $M_2$ with $f(M_1)\subset M_2$. Then
$f$ is locally finite  and
the normal component of $f$ has non vanishing normal derivative
at each point of $M_1$.
 So, by the standard result,
it follows that $f$ is locally proper at each point in $M_1$ and
preserves the sides.

\noindent (7.h) (Bishop's Lemma  [Chi]): Let
$E$ be a (complete) 
pluripolar set in the domain $U\times U'\subset {\bf C}^4$,
where $U$ and $U'$ are domains in $\bf C^2$. Let $A\subset
(U\times U')\sm E$ be 
an analytic set of pure dimension $2$ with no limit points 
in $U\times \PP U'$. Suppose that there exists 
an open subset $V\subset U$ such
 that   $\overline {A}\cap((V\times U'))$ is an analytic variety
in $V\times U'$ of pure dimension 2. 
Then   $\overline {A}$ is an analytic set
in $U\times U'$.

\bigbreak

{\bf Remark} (7.1)$'$(a): We observe that (7.d) 
can be strengthened as follows: For 
any point $\wt{p}(\in T_2^+)\approx p$,
 it also holds that $Q_{\wt p}\cap {\cal U}_2(p)\cap\PP D\subset
 \PP D^+\cup {\wt p}$. 
 The proof of this fact is the same as that for (7.d) (see [DP]).
 However, since this fact will be important in our later discussion,
 for the convenience of a reader, we give the following details (see 
 a similar argument in [DP]):
  
After an appropriate holomorphic change of coordinates,
 we can assume that $p=0$ and $T_2^+$ near
 $p$ is $i{\bf R}^2=\{(iy_1,iy_2)\}$. We also choose $z_1$-axis
for the complex tangent space of $\PP D$ at $0$.
 By the way we chose $T^+_2$, we can assume that $D$ is defined near
 $0$ by $\rho=2x_2+(2x_1)^{k}a(z_1,\overline {z_1},y_2)$, where $m>1$ and
 $a|_{i{\bf R}^2}\not\equiv 0$. By the way we chose $T^+_2$,
we can see that $k$ is even and $a(0)>0$ (see [DP]).
Now, for $w=(ib_1,ib_2)\in i{\bf R}^2$, $Q_{w}$ is defined by
the equation: $z_2-ib_2+(z_1-ib_1)^{2m}a(z_, \overline {w})=0.$ So, 
it can be seen that
near the origin, $Q_{w}\cap i{\bf R}^2=\{w\}$. 
This completes the proof of the
above fact.

(b)
To see $(7.a)$, 
we notice the facts that $Q_z\cap \PP D$
is at most a real analytic subset of dimension 1
 and $T_1\cup T_0$ is contained in
a real analytic set of at most dimension 1 [BM].
Therefore, ${\cal U}_2(p_0)\cap \PP D\cap Q_{T}$ is contained in a real
 analytic
space which can be fibered along a
 one dimensional real analytic curve and each fiber
has at most dimension 1. Thus, 
 ${\cal U}_2(p_0)\cap \PP D\cap Q_{T}$ has at most dimension 2.

 (c) (7.e) can be seen as follows: For each $\wt p$ close to $p$, after
 a normal
holomorphic change of coordinates (see related discussions in $\S 2$)
and using the strong pseudoconvexity of $\PP D$ at $\wt p$, 
one sees clearly that $Q_{\wt p}$ cuts $\PP D$ only at $\wt p$. Since this
process can be done uniformly, we can achieve (7.e) after shrinking the size
of ${\cal U}_2$.

\bigbreak  
{\bf Lemma 7.2}: Shrinking ${\cal U}_1(p)$ if necessary,
 one can 
make ${\cal V}(p,q)$ having no limit points 
in $({\cal U}_1(p)\sm {\overline {D}})\times \Sigma'$.
\bigbreak
{\it Proof of Lemma 7.2}: Suppose that $(w_v,w'_v)\in {\cal V}$
with $(w_v,w'_v)\ra (w_0,w_0')\in({\cal U}_1\sm \overline {D})\times \Sigma'$. Then
$f(Q_{w_v}\cap {\cal U}_2\cap D)\supset\ _{w'{}^*_v}Q'_{w'_v}$. 
In particular, from the properness of $f$, it can be seen that
there exists
a point $\xi^*_{v}\in Q_{w_v}\cap \overline {\cal U}_2\cap D$ 
such that $f(\xi^*_v)= w^{'*}_v$ and
$f(Q_{w_v}\cap O(\xi^*_v))\supset _{w'{}^*_v}Q'_{w'_v}$ Assume, 
without loss of generality, that $\xi^*_{v}\ra \xi_0\in \overline {\cal U}_2\sm D$ 
and $w^*_v\ra w'_0$. 
Then $f(\xi_0)=w'_0$ and $\xi_0\in\PP D$. 
By Lemma 7.1, it follows that $f$ extends holomorphically to a neighborhood of 
$\xi_0$. Therefore, there is
a small open subset $\O$ near $\xi_0$ such that 
$f$ is a proper holomorphic map
from $\O$ to $f(\O)$ and ${f|_{\overline {\O}}}^{-1}(f(\xi_0))=\xi_0$
 (see (7.f)). Meanwhile, $f(\O\cap D)\subset D'$,
$f(\O\cap D^c)\subset D'{}^c\cap {\cal U}_2'$. 
Now, let $\eta_v\in
\O$ be such that $f(\eta_v)=w_v'$. Then, for $v\gg 1$, it holds that
$f(Q_{\eta_v}\cap O({\eta_v^*}))\subset Q'_{w_v'}$. Obviously, $\eta_v\ra \xi_0$ and
 $f(Q_{\eta_v}\cap {\O})\supset f(Q_{w_v}\cap {O}({\xi_0}))$. Write
${f|_{\O}}^{-1}=\{\sigma_1,\cdots,\sigma_k\}$. Hence, we see that
$Q_{w_v}\cap O({\xi_0})\subset \cup_{j=1}^{k}\sigma_j(f(Q_{\eta_v}\cap\O))$.

Similar to Claim 5.2,
we have 
the following
\bigbreak
{\bf Claim 7.2$'$}: For $v$ sufficiently large, 
 $Q_{\eta_v}\cap {\cal U}_2=Q_{w_v}\cap {\cal U}_2$. 
Thus, by passing to the limit, it holds that $w_0\in A_{\xi_0}$.
\bigbreak
{\it Proof of Claim 7.2$'$}:
Let $\rho$ and $\rho'$ 
be the defining functions of $D$ and $D'$, respectively.
Then we notice that $\rho'(f(z),\overline f(w))
=\rho(z,\overline w)h(z,\overline w)$ for $z,w\in \O$  (see (7.f)) 
and $h\not =0$ (we might need to shrink $\O$ here). 

For $z \in\cup_{j=1}^{k}\sigma_j(f(Q_{\eta_v}\cap\O))$, there is a point
 $\wt{z}
 \in Q_{\eta_v}\cap\O$ such that $f(z)=f(\wt{z})$. Now,
 from $$\rho'(f(z),\overline {f(\eta_v)}))=\rho(z,\overline {\eta_v})h(z,\overline {\eta_v})$$
 and 
  $\rho'(f(\wt{z}),\overline {f(\eta_v)})=\rho(\wt{z},\overline {\eta_v})h(\wt{z},
\overline {\eta_v})=0$,
  it follows that $\rho(z,\overline {\eta_v})=0$. Thus, we conclude that
  $z\in Q_{\eta_v}$. Since both $Q_{w_v}$ and $Q_{\eta_v}$ are
  connected complex curves in  ${\cal U}_2$,
  we see the 
$Q_{\eta_v}\cap {O}({\xi_0})=Q_{w_v}\cap {O}({\xi_0})$

Passing to a limit, we now conclude that
$Q_{\xi_0}\cap {\cal U}_2=Q_{w_0}\cap {\cal U}_2$.
\bigbreak
Next,
$w_0\in A_{\xi_0}$. Since $\xi_0\in\PP D$,  $w_0\in\PP D$. 
This is a contradiction and thus completes the proof
of Lemma 7.2. $\endpf$
\bigbreak
In what follows, we always make ${\cal U}_1$ small
 such that Lemma 7.2 holds.
Following the strategy of Diederich-Pinchuk,
we now prove the following lemma.
\bigbreak
{\bf Lemma 7.3}: Suppose that for certain  ${\cal U}_2(p)$  small enough, 
it holds that $(Q_{p}\cap \PP D\cap {\cal U}_2(p))\sm
 \{p\}\subset\Sigma$. 
Then after shrinking ${\cal U}_1$ and ${\cal U}'_1$ if necessary,
${\cal V}(p,q)$ is an analytic sub-variety of 
$({\cal U}_1(p)\sm \overline {D})\times ({\cal U}'_1(p) \sm \overline {D'})$
which has no limit points in
$\left(({\cal U}_1\sm \overline {D})\times \PP {\cal U}_1'\right)\cup
\left( ({\cal U}_1\sm \overline {D})
\times \Sigma'\right)$.
 We mention that in this lemma,
$(Q_{p}\cap \PP D\cap {\cal U}_2(p))\sm \{p\}$ can be empty, and
all Segre surfaces are defined in the Segre system
as chosen here.
\bigbreak
{\it Proof of Lemma 7.3}: Part of the proof 
is greatly motivated by an argument in [DP] for the study of biholomorphic
mappings.

   Let ${\cal U}_2$ be sufficiently small  so that $(Q_{p}\cap 
   \overline { D}\cap \overline {\cal U}_2)\sm \{p\}\subset\Sigma$. We also
assume the property in (7.c).
   When
$(Q_{p}\cap 
   \overline { D}\cap \PP{\cal U}_2)\sm \{p\}\not =\emptyset$, 
then by the hypothesis,
there is a small $\epsilon >0$
such that $f$ extends holomorphically
to $\overline {B_{3\epsilon}(E_p)}$, where
$E_p=Q_p\cap \overline {D}\cap \PP {\cal U}_2$ and
$B_{\delta}(E_p)=\{z:\ \h{dist}(z,E_p)<\delta\}$ for $\delta>0$.
Moreover, by (7.f), we can assume that $f$ is  finite to one from
$B_{3\epsilon}(E_p)$ and $p\not\in B_{4\epsilon}(E_p)$.
Now, we
let $R'(p)=\{w'\in 
{\cal U}_1'\sm { D'}:\ f(Q_{p}\cap B_{2\epsilon}(E_p))\supset
\ _{w'{}^*}Q'_{w'}\}$. Then we claim that $R'(p)$ is a finite set.
 Indeed, this can be seen as follows: Let 
$\{w'_v(\in R'(p))\}$ with $w'{}^*_v\ra w'{}^*$. Choose $w_v\in  
Q_{p}\cap B_{2\epsilon}(E_p)$ such that $f(w_v)= {w'}_v^*$.
 Assume also 
 that $w_v\ra w\in Q_{p}\cap \overline {B_{2\epsilon}(E_p)}$. 
Then as above,
$f$ is a proper
 holomorphic map from a small
 neighborhood of $w$.
 This implies that
 $f(Q_{p}\cap O(w))$ is an 
 analytic curve near $w'{}^*$. Thus, it follows that 
 any $Q'_{w'_v}$ coincides with $Q'_{w'}$ when $v\gg 1$.
 So, by the 
finiteness of $A'_{w'}$, it follows that $\{w'_v\}$ is a finite sequence.  

Let
$R(p)=\{w\in Q_{p}\cap B_{2\epsilon}(E_p): \ f(w)=w'{}^*\ \h{for some}\ 
w'\in R'(p)\}$. Then $R(p)$ is also finite.
 Slightly shrinking ${\cal U}_2$ and $\epsilon$
 if
necessary, we can assume that 
$R(p)\cap \PP {\cal U}_2\cap \overline {D}=\emptyset$.

To complete the proof of Lemma 7.3, we
need to  show that after shrinking ${\cal U}_1$,
${\cal V}$ is an analytic sub-variety in 
$({\cal U}_1\sm \overline {D})\times
({\cal U}'_1\sm \overline {D'})$, which has no limit points in 
$\left(({\cal U}_1\sm \overline {D})\times \PP {\cal U}_1'\right)\cup
\left( ({\cal U}_1\sm \overline {D})
\times \Sigma'\right)$.

We  first claim that $\cal V$ is a local variety.
To  this aim,
we will prove  that after shrinking the size of ${\cal U}_1$,
for any $(w,w')\in {\cal V}$,
then $E^*(w,w')\cap \PP{\cal U}_2\cap Q_w \cap \overline {D}=\emptyset$,
where
 $E^*(w,w')=\{\xi\in ({\cal U}_2\cap D)\cup B_{\epsilon}(E_w)\cap Q_{w}:
f( _{\xi}Q_{w})\supset _{w^{'*}}Q'_{w'}\}$. 
Here $E_{w}=Q_w\cap \overline {D}\cap\PP {\cal U}_2$ and $B_{\epsilon}(E_{w})$
is defined in a similar way.
Also, by $f( _{\xi}Q_{w})\supset _{w^{'*}}Q'_{w'}$, we mean that for
each small neighborhood $O(\xi)$ of $\xi$, it holds that
 $f(Q_{w}\cap O(\xi))\supset _{w^{'*}}Q'_{w'}$.

 We also make $\overline {B_{\epsilon}(E_w)}\subset
B_{2\epsilon}(E_p)$  by shrinking ${\cal U}_1$.

Indeed,
   when
$E_p=(Q_{p}\cap 
   \overline { D}\cap \PP{\cal U}_2) =\emptyset$, a continuity argument then shows
that
$E_w=(Q_{w}\cap 
   \overline { D}\cap \PP{\cal U}_2) =\emptyset$ for $w\approx p$.
 
Now, if the above assertion does not hold
in case 
 $E_p\not =\emptyset$, after shrinking ${\cal U}_1$ again,
we can then find several sequences
$\{z_j\}\subset \PP {\cal U}_2\cap \overline {D}\cap E^*(w_j,w'_j)$
 with $z_j\ra z_0$,
$\{w_j\}\subset {\cal U}_1\sm \overline {D}$ with $w_j\ra p$, and
$\{w_j'\}\subset {\cal U}_1'\sm\overline {D}'$ with $w_j'\ra w_0'\in {\cal U}_1'$
 such that the following holds:

\indent (a) $z_j\in Q_{w_j}$,
(b) $f(z_j)=w_j'$,  and (c) $f(_{z_j}Q_{w_j}\cap B_{\epsilon}(E_{w_j}))\supset 
_{w_j^{'*}}Q_{w_j'}$ when $j>>1$.
We observe that  $f$ extends  to a proper holomorphic map
from $O(z_0)\cap B_{\epsilon}(E_p)$,
Hence, by passing to the limit,
 we can see that
$f( _{z_0}Q_p\cap B_{\epsilon}(E_p))\supset _{w_0^{'*}}Q'_{w'_0}$,
$f(z_0)=w'_0$, and
$z_0\in Q_p\cap\PP {\cal U}_2\cap\overline {D}$. Hence, we have a contradiction;
for it implies that $z_0\in R(p)\cap \PP {\cal U}_2$.

Now, by making the size of ${\cal U}_1$ sufficiently small,
 we can assume that for each $(w,w')\in {\cal V}$,
 $E^*(w,w')\cap D\cap \overline {\cal U}_2\subset {\cal U}_2\cap D$.
Write  $E^*(w,w')\cap D\cap \overline {\cal U}_2=\{a_1,\cdots,a_m\}\subset D\cap
 {\cal U}_2$.
Let $(\wt w, \wt w')
(\in {\cal V})\approx (w,w')$. Then, since $E(\wt{w},\wt{w}')\cap D\cap \PP
{\cal U}_2=\emptyset$, it can be seen that
there exists  a point 
$\wt \xi$
 so that $\wt {\xi}\in Q_{\wt w}\cap {\cal U}_2\cap D$,
  $f(_{\wt \xi}Q_{\wt w})
\supset _{\wt{w'}^*}Q'_{w'}$, and $f(\wt \xi)=\wt{w'}^*$. Indeed,
let $\{\wt \xi_j\}=f^{-1}(\wt {w^{'*}})\cap \overline {\cal U}_2\cap D$.
Then since $f$ is locally proper near each $\wt{\xi_j}$, the choice
of our $w'$ indicates that $f(O(\wt{\xi_l})\cap Q_{\wt w})\supset
_{\wt{w'}^*}Q'_{w'}$ for a certain $l$. So, we can take
$\wt \xi$ to be this $\wt \xi_l$, which has to stay in ${\cal U}_2\cap D$.

By the properness of $f$ and
by the continuous dependence of Segre surfaces on the base points,
it follows   that
$\wt{\xi}$ has to be close to one of the $a_j's$ when $\wt w$ is close 
to $w$; for otherwise we would have
$E^*(w,w')\cap \overline {D}\cap\PP {\cal U}_2\not =\emptyset$. 
We let $\rho_2$ be a real analytic defining function of
$\PP D'$ near $q=f(p)$, and we  choose
a coordinates system near $p$ so that the tangent space $T^{(1,0)}_{p}
\PP D$ is the $z_1$-axis. Then $Q_{\wt w}$ is given by 
$z_2=h(z_1,\overline {\wt w})$ for some holomorphic function $h$. We claim that
the  condition
that $f(Q_{\wt w}\cap D\cap {\cal U}_2)\supset \ _{\wt{w'}^*}Q'_{\wt w'}$
for $(\wt{w},\wt{w'})\approx (w,w')$,
 can be expressed by the equation $\rho_2(f(z),\overline {\wt w}')=0$ with
$z\in Q_{\wt w}\cap O(a_j)$ for some $j$.  Indeed, 
 $\rho_2
(f(z),\overline {\wt w}')=0$ implies that $f(Q_{\wt w}\cap O(a_j))\subset Q'_{\wt w'}$.
Since we can make $O(a_j)$  sufficiently small so that $f$ is
proper from $O(a_j)$ to a neighborhood, say $O_j'$, of $w^{'*}$,
 $f(Q_z\cap O(a_j))$ is also an analytic variety of $O_j'$ of codimension 1.
 Since we can choose  $O_j$ and $O_j'$ suitably so that $O_j'\cap Q'_{\wt w'}$ is connected
 for $\wt{w'}\approx w'$, 
thus  $f(Q_{\wt w}\cap O(a_j))=Q_{\wt w'}\cap O_j'$.
 Now, if we make $\wt{w'}$ close
enough to $w^{'}$, then by the continuity of the ${\cal R}$-operator,
 $\wt{w}^{'*}\in  f(Q_{\wt w}\cap O(a_j))$. This verifies our
claim.

On the other hand, we have the following Taylor expansion near $a_l^{(1)}$
for each $l$, where
$a_l^{(1)}$ is the first coordinate   of $a_l$:
$$\rho_2(f(z),\overline {\wt w}')=\sum_jb^l_j(\overline {\wt w},\overline {\wt w'})(z_1-a_l^{(1)})^j$$
Let $V_l$ be the common zeros of $b_j^{l}(\overline {\wt w},\overline {\wt w'})$ for
$j=0,1,\cdots$. Then clearly near $(w,w')$, ${\cal V}$ is the finite union
of $V_l$'s. This shows that ${\cal V}$ is a local variety.

We next let $(w_v,w'_v)(\in {\cal V})\ra
(w_0,w'_0)\in ({\cal U}_1\sm \overline D)\times ({\cal U}_1'\sm \overline {D}')$. 
 We then  show that
$(w_0,w_0')\in {\cal V}$, too.  Notice that
$f(Q_{w_v}\cap {\cal U}_2)\supset\ _{{w'}_v^*}Q'_{w_v'}$. So, as argued before,
 there is a point 
$\xi_v\in Q_{w_v}\cap {\cal U}_2\cap D$ such that $f( _{\xi_{v}}Q_{w_v})
\supset_{w^{'*}_v}Q'_{w'_v}$.
Assume, without loss of generality, that $\xi_0=\lim \xi_v$.
Letting $v\ra \infty$, since $\xi_0\in Q_{w_0}\cap \overline {\cal U}_2\cap {D}$ and since $f$ is locally proper near $\xi_0$,
it follows that
 $f(O(\xi_0)\cap Q_{w_0}\cap (B_{\epsilon}(E_{w_0})\cup ({\cal U}_2\cap D))
\supset _{w_0^{'*}}Q'_{w_0'}$.
 Thus, by the above choice of  
${\cal U}_2$, it follows that $\xi_0\in {\cal U}_2$. We therefore
 conclude that
$f(Q_{w_0}\cap {\cal U}_2\cap D)\supset\ _{{w'}_0^*}Q'_{w_0'}$.
This shows  that $(w_0,w_0')\in 
 ({\cal U}_1\sm \overline {D})\times ({\cal U}'_1\sm \overline {D'})$.

 Next,  applying Lemma 7.2 and (7.b), we 
 conclude  that ${\cal V}$ has no limit points in 
 $  ({\cal U}_1\sm \overline {D})\times \PP {\cal U}'_1\cup
 ({\cal U}_1\sm \overline {D})\times \Sigma'$.  $\endpf$

\bigbreak
We still keep the above notation and assume the hypothesis in Lemma 7.3.
Notice that $\cal V$ contains points of dimension 2 
(see the observation made after Lemma 7.1).
 We see that $ {\cal V}$ is an analytic variety of dimension 2 in 
  $\left(({\cal U}_1\sm {\overline  D})\times {\cal
 U}'_1\right)\sm {A})$. Here $A=({\cal U}_1\sm \overline {D})
 \times (T'_1\cup T'_0)$.
 
 Decompose $\cal V$ into a locally finite union of its 
irreducible components. In what follows, we will be only
 interested in the locally finite union of those two dimensional components,
which we will still denote by $\cal V$, for the brevity of notations.
(The new $\cal V$ can be formed by considering
 all smooth points of dimension 2 in the old $\cal V$ and then adding 
new points by taking the
 closure). Now,
we can say that ${\cal V}$ is an analytic variety in
$(({\cal U}_1\sm \overline {D})\times {\cal U}'_1))\sm {A}$,
which has pure dimension 2.

Next, we let 
 $\wt{A}=({\cal U}_1\sm \overline {D})\times \wt{T}(\supset A)$,
which
 is a  pluripolar set in $({\cal U}_1\sm \overline {D})\times {\cal U}'_1$ by (7.b).
 We also  write 
$\wt{\cal V}={\cal V}\sm \left(({\cal U}_1\sm \overline {D})\times \wt{T})\right)$.
We now  think of ${\cal V}$ and $\wt{\cal V}$ as  
subsets in $({\cal U}_1\sm {\overline  D})\times 
{\cal
 U}'_1$. 
 Then, 
 we see that $\wt {\cal V}$ is an analytic variety of pure
 dimension 2 in 
  $\left(({\cal U}_1\sm {\overline  D})\times {\cal
 U}'_1\right)\sm \wt {A}$. 

\bigbreak
{\bf Lemma 7.4}: 
As before, let $\{{\cal U}_j\}_{j=1}^{3}$ be a Segre
neighborhood system at $p$ and 
 assume  the hypothesis in Lemma 7.3.
If there exists a sequence
$\{p_j\}\subset \Sigma$ with $p_j\ra p$,
 such that $Q_{p_j}\cap \PP D\cap {\cal U}_2\subset \Sigma$, then
${\cal V}$ extends as an analytic variety of pure dimension $2$ in
${\cal U}_1\times {\cal U}'_1$, after shrinking
${\cal U}_1$ and ${\cal U}'_1$. Thus, by Theorem 1.1, $f$ admits
a holomorphic extension at $p$.
\bigbreak
{\it Proof of Lemma 7.4}: The main idea of the proof is to apply the 
Bishop lemma, as first used in the work of Diederich-Pinchuk  [DP].
 Without loss of generality, we
also assume that any CR function defined near $p\in \PP D$ can  be
holomorphically extended to $D$ (near $p$) (see [BT] and [Tr]).  

We let $ z$ be a certain $p_j$ with $j\gg 1$. Slightly shrinking
${\cal U}_2$, we observe that 
 $Q_{ z}\cap \PP D\cap \overline 
{\cal U}_2\subset \Sigma$. It thus
follows that
for a small neighborhood $G$ of $z$, one also has
 $Q_{\wt z}\cap \PP D\cap \overline {\cal U}_2\subset \Sigma$ for any $\wt z\in\overline { G}$.
 Now, we consider ${\cal V}^*
 =\left({\cal V}\cap ((G\sm \overline {D})\times {\cal U}'_1)\right)$. 
 Then as did in Lemma 7.2, one sees that ${\cal V}^*$ has no limit
points in $(G\sm \overline {D})\times \Sigma'$. 
We claim that ${\cal V}^*$ has no limit points
in $(G\sm \overline {D})\times (T'_1\cup T'_0)$, neither. 
To see this, we let $(w_v,w_v')
(\in {\cal V}^*)\ra (w_0,w_0')\in (G\sm \overline {D})\times (T'_1\cup T'_0)$. 
Then, as in
 Claim 7.2 $'$, one sees that for certain $\xi_0\in Q_{w_0}\cap \PP D\cap \overline {\cal U}_2$, it holds that $f(\xi_0)=w_0'$. 
However, by our choice of $G$, we know that
$\xi_0\in\Sigma$. Thus, a similar argument as in Claim 
7.2$'$ indicates that
$w_0\in A_{\xi_0}\in\PP D$. This is a contradiction.

Hence, from the above argument, it follows that 
$\overline {{\cal V}^*}\cap ((G\sm \overline {D})\times {\cal U}'_1)={\cal V}^*$ is 
an analytic variety of pure
 dimension 2 in $(G\sm \overline {D})\times {\cal U}'_1$ without limit
points in ${\cal U}_1\times \PP {\cal U}'_1$. 
Now, we let $\wt{\cal V}^*={\cal V}^*\sm \wt{A}$. Then
$\overline {\wt{\cal V}^*}\cap((G\sm \overline {D})\times{\cal U}'_1)
=
\overline {\wt{\cal V}}\cap((G\sm \overline {D})\times{\cal U}'_1)$. 
\bigbreak
{\bf Claim 7.5}: $\overline {\wt{\cal V}^*}\cap((G\sm \overline {D})\times{\cal U}'_1) 
={\cal V}^*$. Thus,
$\overline {\wt{\cal V}}\cap((G\sm \overline {D})\times{\cal U}'_1) $ is an analytic
variety of pure dimension 2 in $(G\sm\overline {D})\times {\cal U}'_1$.
\bigbreak
{\it Proof of Claim 7.5}: Let 
$\pi_{G}: {\cal V}^*\ra G$ and $\pi'_G:{\cal V}^*\ra {\cal U}'_1$.
We first show that $\pi'_G$ is finite to one, after shrinking $G$ if
necessary. To this aim, assume that $w'\in {\cal U}'_1$ and
let $\{z_j\}$ be a sequence in $\pi_G({\pi'_G}^{-1}(w'))$ 
with $z_j\ra z_0(\in \overline {G})$.
As in Lemma 7.3, we can then assume that for a
certain sequence $\{\xi_j(\in Q_{z_j})\}$, 
it holds that
$\xi_j(\in Q_{z_j}\cap\overline {\cal U}_2)\ra \xi_0\in Q_{z_0}\cap \overline {\cal U}_2\cap\overline {D}$, 
$f(\xi_j)\ra w'{}^*$, and $f(_{\xi_j}Q_{z_j})
\supset \ _{{w'}^*}Q'_{w'}$. As before,
by our choice of $G$, if follows that there is a small neighborhood
of $\xi_0$ which is mapped properly by $f$ to a neighborhood of $w'{}^*$.
Thus, $f^{-1}(\ _{w'{}^*}Q_{w'})$ has only finitely many components near 
$\xi_0$. Therefore, we conclude that $\{Q_{z_j}\}$ is a finite sequence.
Making use of the fact that $A_{z}$ is finite for each $z\in G$, we
see that $\pi_{G}^{-1}(w')$ is a finite set and
thus $\pi_{G}$ is a local analytic covering map. 
In particular, $\pi'_{G}$ is  an open mapping.

Next, we consider
${{\cal V}^*}\cap((G\sm \overline {D})\times \wt{T}) $. Then the above 
fact indicates that it does not contain any non-empty open subset
of $ {\cal V}^*$.
In fact, if this is not the case,  $\pi'_{G}(
{{\cal V}^*}\cap((G\sm \overline {D})\times\wt{T})\subset \wt{T} $ 
would contain open subsets of ${\cal U}'$, which contradicts
the property of $\wt T$ (see (7.b)). 
Hence, we see that
 $\overline {\wt{\cal V}^*}\cap((G\sm \overline {D})\times{\cal U}'_1) 
={\cal V}^*$. $\endpf$
\bigbreak

Now, we can apply the Bishop extension lemma
(7.h) to conclude that $\wt{\cal V}$ extends as an analytic 
variety of dimension
2 in $({\cal U}_1\sm \overline {D})\times {\cal U}'_1$, 
 which has no limit points
in ${\cal U}_1\times \PP {\cal U}'_1$. Let us 
 still write the extended variety as ${\cal V}$.
 Since any compact analytic variety in $\bf C^2$ has to be a finite set,
 we can see, by using the above mentioned property of $\cal V$, 
 that
  $\pi^{-1}(z)$  is finite for
 each $z\in {\cal U}_1\sm\overline {D}$, and 
$\pi$ is proper. Hence, ${\cal V}$ can have only finitely many 
 irreducible components, $\{{\cal V}\}_{j=1}^m$, (of dimension 2); and
 each  of them has to be a sheeted analytic covering space
over ${\cal U}_1\sm \overline {D}$ ([Wh]).

 So, each ${\cal V}_i$ can be defined by equations of the form
([Wh]):
$$P_{ik}(w,w'_k)=w_k'{}^{N_{ik}}+\sum_{0\le j<N_{ik}}a_{ikj}(w){w'}_k^{j}$$
 with $a_{ikj}(w)$
 holomorphic and bounded in ${\cal U}_1\sm \overline {D}$. By the result
 of Baouendi-Treves [BT] and Trepreau [Tr], we see that $a_{ikj}$ can be
 holomorphically extended to a neighborhood of $p$. So
each ${\cal V}_i$ extends to a neighborhood of $p$. Notice that
 there are points arbitrarily close to $p$ where $f$ extends, we conclude that $\Gamma_{f}\subset \cup {\cal V}_i$. This tells that $f$ extends as
a holomorphic correspondence
 near $p$. Applying Theorem 1.1, we see that $f$ extends holomorphically 
 across $p$. $\endpf$
\bigbreak
We are now ready to complete the proof of Theorem 1.3$'$ by
the argument appeared in [DFY]:
 First, $\Sigma\supset \PP D^-\cup T_2^-\cup T^{\pm}_2$.
 Since $f$ extends almost everywhere, from  (7.e), Lemma 7.3, and Lemma 7.4, it follows that $\PP D^+\subset \Sigma$. Applying Lemma 7.1,
 Lemma 7.3 and Lemma 7.4 with
Remark 7.1$'$(a), we conclude, 
from the fact $\PP D^+\subset \Sigma$, that $T^+_2\subset \Sigma$. Now,
 we can use (7.a) to show that
$T_1\cup T_0\subset\Sigma$:

 Let $z (\approx p_0)\in T_1$.
Since $A_{p_0}$ is finite, it follows that only finitely many
$z$'s can make $Q_{z}\cap \PP D$ to be a real analytic set of dimension 1
 in $T_1$. Thus, by
Lemma 7.3, Lemma 7.4, and (7.a); it follows that $f$ extends at most except at finitely many points near $p_0$.
Applying again Lemma 7.3 and Lemma 7.4, we see that $f$ extends also
 as a holomorphic correspondence at those
isolated points, too.  The proof of Theorem 1.3$'$ 
is now complete. $\endpf$
\bigbreak
 {\bf Remark 7.6}: Except the first part of Lemma 7.1, the entire
proof of Theorem 1.3$'$
 is purely a local argument. So, using the following Lemma $7.1''$ to
replace  part of Lemma 7.1, we also see the proof of the following
local version of Theorem 1.3
\bigbreak
{\bf Theorem 1.3}$''$: Let $D_1$ and $D_2$ be two bounded domains in $\bf C^2$.
Assume that $M_1$ and $M_2$ are contained in the real analytic 
boundaries of $D_1$ and $D_2$, respectively.
 Suppose that $f$ is
a proper holomorphic mapping from $D_1$ to $D_2$, that is continuous
on $D_1\cup M_1$ and sends $M_1$ into $M_2$. Also assume that $M_1$ and $M_2$
are of finite type. Then $f$ admits a holomorphic extension across $M_1$.
\bigbreak
{\bf Lemma 7.1}$''$: Assume the notation and hypothesis in Theorem $1.3''$.
Let $p_0\in M_1$ be such that $D_2$ is pseudoconvex  near
$q_0=f(p_0)$. Then $f$ admits a holomorphic extension at $p_0$. 
\bigbreak
{\it Proof of Lemma 7.1$''$}: Construct a pseudoconvex domain of
finite type $\O_2\subset
 D_2$ with $\PP \O_2$ containing a small piece of $M_2$ 
near $q_0$. Let $\O_1$ be a connected component of $f^{-1}(\O_2)$ that
contains $p_0$ in   part of its smooth boundary. Then $f$ is proper from
$\O_1$ to $\O_2$ and maps $M_1$ near $p_0$ to $M_2$ near $q_0$. Meanwhile,  
it also follows that $\O_1$ is pseudoconvex.
Now, by
examining the proof of Theorem 2 of [BC] and by [BBR], to see the extension
of $f$ at $p_0$, it suffices for us to show the for $z (\approx z_0)
\in \O_1$,
it holds that $\h{dist}(f(z),M_2)^\g\ale \h{dist}(z,M_1)$, for some positive
integer $\g$. But this can be seen by applying the Hopf lemma, 
a classical argument 
of Henkin,
and using the existence of family of plurisubharmonic
peaking functions near $p_0$; as did  in  [Ber] (pp 629, line -15, pp 629, 
line 11) and [BC] (pp359  line -6, pp 360 line 7).

\end{document}